\documentclass[11pt]{amsart}

\usepackage{amsfonts}
\usepackage{amsmath}
\usepackage{amsthm}
\usepackage{mathrsfs}
\usepackage{enumerate}
\usepackage{graphicx}
\usepackage{tikz}
\usepackage{caption}
\usepackage{subcaption}
\usepackage{caption}
\usepackage[square,sort,comma,numbers]{natbib}
\usepackage{color}
\usepackage{hyperref}
\usepackage{setspace}
\usepackage{natbib}
\usepackage[T1]{fontenc}
\usepackage{times}
\usepackage[left=2cm,right=2cm,top=2.5cm,bottom=2.5cm]{geometry}

\usepackage{bigints}

\newcommand{\EXCLUDE}[1]{}

\newcommand{\remove}[1]{}
\newcommand{\beq}{\begin{eqnarray}}
	\newcommand{\eeq}{\end{eqnarray}}
\newcommand{\beqq}{\begin{eqnarray*}}
	\newcommand{\eeqq}{\end{eqnarray*}}%\pagestyle{plain}
\def\:{:\,}
\DeclareMathOperator*{\argmax}{arg\,max}

\newtheorem{thm}{Theorem}[section]
\newtheorem{cor}[thm]{Corollary}
\newtheorem{lem}[thm]{Lemma}
\newtheorem{prop}[thm]{Proposition}

\numberwithin{equation}{section}

\def\rar{\rightarrow}

\newcommand{\be}{\beta}

\newcommand{\ep}{\epsilon}
\newcommand{\al}{\alpha}

\newcommand{\lam}{\lambda}
\newcommand{\Lam}{\Lambda}

\def\th{\theta}

\def\bR{\mathbb{R}}

\def\bZ{\mathbb{Z}}

\def\bN{\mathbb{N}}

\newcommand{\cA}{{\mathcal A}}
\newcommand{\cB}{{\mathcal B}}
\newcommand{\cE}{{\mathcal E}}
\newcommand{\cP}{\mathcal{P}}

\newcommand{\cH}{{\mathcal H}}

\newcommand{\cC}{{\mathcal C}}

\newcommand{\cT}{{\mathcal T}}

\def\1{\mathbf{1}}

\def\lb{\left(}
\def\rb{\right)}

\def\lab{\label}

\def\nn{\nonumber}
\def\f{\frac}
\def\llra{\longleftrightarrow}

\begin{document}
	
	%%%%%%%%%%%%%%%%%%%%%%%%%%%%%%%%%%%%FRONTPAGE%%%%%%%%%%%%%%%%%%%%%%%%%%%%%%%%%%%%%
	%\bibliographystyle{abbrv}
	
	% "Title of the paper"
	\title[Phase transition and percolation at criticality]{\large{Phase Transitions and Percolation at Criticality in enhanced random connection models}}
\author{Srikanth K. Iyer}
\address{(SKI) Department of Mathematics\\
		Indian Institute of Science\\
		Bangalore, India, Email : \textit{srikiyer@gmail.com}}

\author{Sanjoy Kr. Jhawar}
\address{(SKJ) Department of Mathematics\\
	Indian Institute of Science \\
	Bangalore, India, Email : \textit{sanjayjhawar@iisc.ac.in}}

\thanks{SKI's research was supported in part from Matrics grant from SERB and DST-CAS, SKJ's research was supported by DST-INSPIRE Fellowship. Corresponding author email: sanjayjhawar@iisc.ac.in. }
\begin{abstract}
We study phase transition and percolation at criticality for three random graph models on the plane, viz., the homogeneous and inhomogeneous enhanced random connection models (RCM) and the Poisson stick model. These models are built on a homogeneous Poisson point process $\cP_{\lam}$ in $\bR^2$ of intensity $\lam$. In the homogenous RCM, the vertices at $x,y$ are connected with probability $g(|x-y|)$,  independent of everything else, where $g:[0,\infty) \to [0,1]$ and $| \cdot |$ is the Euclidean norm. In the inhomogenous version of the model, points of $\cP_{\lam}$  are endowed with weights that are non-negative independent random variables with distribution $P(W>w)= w^{-\beta}1_{[1,\infty)}(w)$, $\beta>0$. Vertices located at $x,y$ with weights $W_x,W_y$ are connected with probability $1 - \exp\lb -  \frac{\eta W_xW_y}{|x-y|^{\al}} \rb$, $\eta, \al > 0$, independent of all else. The graphs are enhanced by considering the edges of the graph as straight line segments starting and ending at points of $\cP_{\lam}$. A path in the graph is a continuous curve that is a subset of the union of all these line segments. The Poisson stick model consists of line segments of independent random lengths and orientation with the mid point of each segment located at a distinct point of $\cP_{\lam}$. Intersecting lines form a path in the graph. A graph is said to percolate if there is an infinite connected component or path. We derive conditions for the existence of a phase transition and show that there is no percolation at criticality. 

\end{abstract}
%
%\date{\today}
\maketitle
\noindent\textit {Key words and phrases.} Random Geometric Graphs, Random Connection Model, Enhanced Random Connection Model, Percolation, Phase transition, Continuity of the percolation function.

\noindent\textit{AMS 2010 Subject Classifications.} Primary: 82B43,\,%percolation, 
60D05.\, % Geometric probability and stochastic geometry
Secondary: 05C80 % Random graphs.

\section{Introduction and Main Results}
The study of random graphs started with the pioneering work by Erd\"{o}s and R\'{e}yni \citep{Erdos59}, \citep{Erdos60} and Gilbert \citep{Gilbert59} on the Erd\"{o}s-R\'{e}yni model. The random graph in the Erd\"{o}s-R\'{e}yni  model is constructed on a set of $n$ vertices, for some $n\in \bN$ with an edge drawn between any two pairs of nodes independently with probability $p\in[0,1]$. Detailed work on the Erd\"{o}s-R\'{e}yni graph can be found in \citep{Bollobas2001}, \citep{Hofstad2017}, \citep{Janson2000}. The Bernoulli lattice percolation model on $\bZ^d$ is an extensively studied random graph model \citep{Kesten1980}, \citep{Kesten1982}, \citep{Grimmett1999} where the geometry of the underlying space plays an important role. The vertex set is $\bZ^d$ with an edge between points at Euclidean distance one with probability $p$ independent of other edges. The above geometric model was extended to the continuum by considering a point process in $\bR^d$ with edges between points that are within a Euclidean distance $r > 0$. Such a model has found wide application in modeling ad-hoc wireless networks and sensor networks. This model is called the Boolean model \citep{Gilbert61} or the random geometric graph (RGG). The questions of interest in such applications are percolation, connectivity and coverage, details of which can be found in \citep{Franceschetti2007} and \citep{Haenggi2012}. Rigorous theoretical analysis of the percolation problem in such graphs can be found in \citep{Meester1996} while the monograph \citep{Penrose2003} carries a detailed compilation of the important results on the topic in the sparse, thermodynamic and connectivity regimes. When each point of the underlying point process has an independent subset of $\bR^d$ associated with it, then the union of all such sets is what forms the germ-grain model which is of much interest in stochastic geometry \citep{SchneiderWeil2008}.

The main goal of this paper is to derive conditions under which three network models on the plane exhibit a phase transition and show that under some additional conditions the percolation function is continuous. We shall now describe these models and the problem of interest and discuss some applications. All these three models are constructed over a homogeneous Poisson point process denoted $\cP_{\lam}$ in $\bR^2$ with intensity parameter $\lam>0$. A phase transition refers to the abrupt emergence of an infinite component in the graph, in which case we say that the graph percolates. A phase transition is said to occur if there exists a critical value $\lam_c \in (0, \infty)$ of $\lam$ such that for $\lam > \lam_c$ the random graph under consideration percolates and for $\lam < \lam_c$ the random graph does not percolate. It can be shown using ergodicity that for $\lam > \lam_c$, there is an infinite component with probability one. In many of percolation models in $\bR^d$ it can also be shown that there is an unique infinite component by adapting the Burton-Keane argument  \citep{Grimmett1999}, \citep{Meester1996}. The percolation function refers to the probability that a typical vertex in the graph is part of the infinite component. Percolation is equivalent to the percolation function being positive. The continuity of the percolation function is a problem of much interest in the random graph literature. See for instance, \citep{Duminil-Copin2016} for a new proof showing absence of percolation at criticality in the Bernoulli bond percolation model on $\bZ^d$ with $d=2$. The problem remain open, for instance, for $d = 3$.

The random connection model (RCM) is a generalization of the RGG and a continuum version of the long range percolation on lattices. This model was introduced in \citep{Penrose1991} and later studied in the context of wireless networks. In such networks communication between nodes depend on the distance between the nodes as well as the interference coming from transmissions from other nodes in the network \citep{Meester1995}, \citep{Meester1996}. In the random connection model we consider, the vertex set will be a homogeneous Poisson point process denoted $\cP_{\lam}$ in $\bR^d$ with intensity parameter $\lam$ for some $\lam>0$. An undirected edge denoted $\{x,y\}$ exists between vertices located at $x,y$ with probability $g(|x-y|)$ independent of everything else, where $g : [0,\infty) \to [0,1]$ is non-increasing. We denote this graph by $G_\lam$. 
%A phase transition is said to occur in $G_\lam$ if there exists a $\lam_c \in (0,\infty)$ such that 
%there exists an infinite connected component ($G_\lam$ percolates) with positive probability only if $\lam > \lam_c$.  
\citep{Penrose1991} showed that a phase transition occurs in $G_\lam$ if and only if the connection function satisfies $0<\int_{\bR^d}g(|x|)\,dx<\infty$. 

The first model that we consider is the enhanced RCM. To define this model consider the RCM described above on the plane $(d = 2)$. If an edge $\{x,y\}$ exists in the RCM then we refer to $x$ and $y$ as direct neighbors. We view each edge as a straight line segment and denote it by $\overline{xy}$. For any two edges $\{x_1,x_2\}$ and $\{x_3,x_4\}$ in the RCM that intersect, we say that the vertices $x_1, x_2$  are indirect neighbors of $x_3,x_4$ and vice versa. We will refer to the resulting graph as the {\it enhanced } random connection model (eRCM) and denote it by $G^e_{\lam}$. It will be more useful to think of the eRCM as enhancing the available paths in the network rather than introducing additional edges as can be seen from the following applications. Intersecting edges along a {\it path} in the RCM allow for switching from one path (in the original graph) to another. The eRCM can be considered as a model for a road or a pipeline network where connections are made locally and intersecting roads or pipelines allow the traffic or the fluid to switch paths. The above could also be used as a model for thin slab of porous media where connections between nodes resemble pipes and crossing of these pipes allow the fluid to flow from one pipe to another. An alternate model for road networks was studied in \citep{AldousShun2010}, \citep{Aldous2014}. The construction of an optimal road network by using the trade-off between a measure of shortness of route and normalized network length for a one parameter family of proximity graphs is studied in \citep{AldousShun2010}. In \citep{Aldous2014} the author introduces  scale invariant spatial networks whose primitives are the routes between points on the plane. The problems of interest are the existence and uniqueness of infinite geodesics, continuity of routes as a function of end points and the number of routes between distant sets on the plane.

In this context it would be more appropriate to consider an inhomogeneous version of the eRCM model where each vertex is endowed with a weight that is indicative of the size, importance of a city or town. In the basic inhomogeneous model we consider, the vertex set is an independent marking of $\cP_{\lam}$ with weight distribution given by $P(W > w) = w^{-\be} 1_{[1, \infty)}(w)$ for some $\be > 0$. Vertices located at $x,y \in \bR^2$ endowed with random weights $W_x,W_y$ are connected by an edge independently with probability
\begin{equation}
g(x,y) = 1-\exp\left(-\f{\eta W_x W_y}{|x-y|^{\al}}\right)
\lab{eq:connection_function}
\end{equation}
where $\eta, \al$ are positive constants. The graph thus obtained is then enhanced in the same manner as described above to obtain the inhomogeous eRCM (ieRCM). This is the second model that we shall study in this paper. We denote the random graph obtained in the inhomogeneous RCM and the enhanced inhomogenous RCM by $H_\lam,\, H^e_\lam$ respectively. Percolation properties for inhomogeneous random connection model with this type of inhomogeniety was studied for long range percolation model on $\bZ^d$ by Deprez, Hazra  and W{\"u}thrich in \citep{Deprez2015} and in the continuum for fixed intensity $\lam$ by Deprez and W{\"u}thrich in \citep{Deprez2018}. Phase transition is expressed in terms of the parameter $\eta$ instead of $\lam$. A simple scaling argument shows that these two are equivalent. In both these models, a phase transition is shown to occur for $d=1$ only if $\al \beta > 2$ and $1 < \al < 2$ and for $d \geq 2$ only if $\al>d$ and $\al\be>2d$.  For all $d$ the percolation function has been shown to be continuous under the condition that $\al\be>2d$ and $\al\in(d,2d)$. For $d \geq 2$ the case when $\min\{\al,\al\be\}>2d$ is open.

The third random graph model on the plane we shall analyze is the Poisson stick model. This is an example of a model that satisfies the axiomatic conditions of so called scale invariant spatial networks mentioned above. This model which was introduced in \citep{RRoy1991} consists of sticks of independent random lengths whose mid points are located at points of $\cP_{\lam}$ with each stick having a random independent orientation. The half-length of the sticks were assumed to have a density with bounded support. Two points in $\cP_{\lam}$ are neighbors in the resulting graph provided the corresponding sticks intersect. A phase transition was shown to occur in such a graph. The Poisson stick graph appears to be a natural model for a network structure formed by silicon nanowires and carbon and other nanotubes on the surface of substrates.  Percolation, conductance and many other significant properties of these nanowire networks are studied in \citep{Pike1974},\citep{Balberg1983},\citep{Oskouyi2014}, \citep{Serre2013}, \citep{Hu2004}. In this paper we consider the Poisson stick model with stick-length distribution having unbounded support and orientation distributed according to some arbitrary non-degenerate distribution. We study existence of phase transition and the continuity of the percolation function. The enhancement in our first two models has similarity with the Steiner tree. The Poisson stick model is similar to the Poisson line process on the plane. The Poisson line process and the Steiner tree have been used by Aldous and Kendall in \citep{AldousKendall2008} to establish asymptotics of excess route length in arbitrary graphs.

\subsection{Notations} We gather much of the notations we need here for easy reference. We define the notations with reference to the RCM and eRCM. However, they carry over to the ieRCM and the Poisson stick models in the obvious way. Let $C(x)$ ($C^e(x)$) be the connected component containing $x\in \cP_\lam$ in $G_\lam$ ($G^e_\lam$). Without loss of generality we assume that there is a vertex at the origin $O$, that is, we consider the process $\cP_\lam$ under the Palm measure $P^o$, the probability distribution conditioned on a point being at origin. The distribution of $\cP_\lam$ under $P^o$ is the same as that of  $\cP_\lam \cup \{O\}$ under the original measure $P$. 
%In what follows we will write $P$ for both the original measure as well as $P^o$. 
Let $C:=C(O)$,\,$C^e:=C^e(O)$ and define the percolation probabilities for $G_\lam,\, G_\lam^e$ as
\begin{equation}
\th(\lam):=P^o(|C|=\infty) \mbox{ and } \th^e(\lam):=P^o(|C^e|=\infty).
\end{equation}
The percolation thresholds denoted by $\lam_c,\lam^e_c$ for the graphs $G_{\lam}$ and $G^e_{\lam}$ respectively  are defined as
\begin{equation}
\lam_c:= \inf\{\lam > 0 : \th(\lam)>0\} \mbox{ and }\lam^e_c:=\inf\{\lam > 0: \th^e(\lam)>0\}.
\end{equation}
Similarly let $\tilde{\lam}_c,\tilde{\lam}^e_c, \lam_{PS}$ be the percolation thresholds for the random graphs $H_{\lam}, H^e_{\lam}$ and $PS_\lam$ respectively.

For any connected region $D\subset \bR^2$ an event $E$ is said to be {\it $D$-measurable}  provided the occurrence or otherwise of $E$ does not depend on the points of $\cP_{\lam}$ that fall outside $D$. 

As mentioned earlier, since we work with paths in the graph, we shall often view the enhanced models as providing additional paths in the original graphs rather than adding edges. In this view, edges in the graphs $G_{\lam} , H_{\lam}$ are straight line segments joining the vertices of $\cP_{\lam}$. Given any of the graphs $G^e_{\lam}, H^e_{\lam}$ or $PS_\lam$ and $x,y \in \bR^2$ we say that there is a {\it path} from $x$ to $y$ if there exists a closed continuous curve from $x$ to $y$ contained entirely in $\cup_{i=1}^{n}e_i$ for some edges $e_1,e_2\cdots,e_n$  in the case of $G_{\lam} , H_{\lam}$ and sticks in the case of $PS_\lam$. Paths thus need not start or end at vertices in the graph.

\begin{figure}
  \centering
  \includegraphics[width=0.5\linewidth]{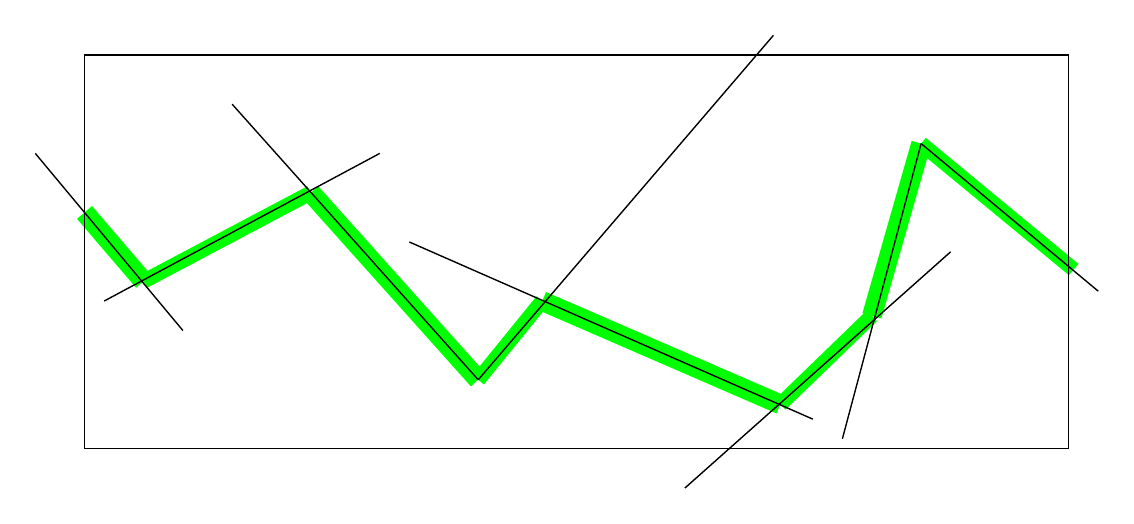}
  \captionof{figure}{The green curve is a left-right crossing of the box.}
  \label{fig:Path}
\end{figure}

A path is said to cross a box $[a,b]\times [c,d]$ if the path is completely contained within the box with end points on opposite sides. We shall refer to these paths as crossings (see Figure~\ref{fig:Path}). 

{\bf Crossing events:} For $s > 0$ and $\rho > 1$ let $LR_s(\rho)$ be the event that there exists a crossing along the longer side of the rectangle $[0,\rho s]\times [0,s]$ and $TD_s(\rho)$ be the event that there exists a crossing along the shorter side of the rectangle $[0,\rho s]\times [0, s]$. $C_s(\rho) :=P(LR_s(\rho)) $.
% which by rotation invariance equals $P(TD_{\rho s}(\f{1}{\rho}))$  \red{not true}

A circuit around $S$ in the region $T\setminus S$, $S\subset T\subset \bR^2$ where both $S,T$ are connected subsets, is a path that starts and ends at the same point and is entirely contained inside $T\setminus S$. Note that the end points of the edges that contribute to the path may lie outside $T\setminus S$. Let $B_s:=[-s,s]^2$ and for $t> s$ let $A_{s,t}:= B_t \setminus B_s$. $\cA_s$ be the event that there exists a circuit in the annulus $A_{s,2s}$. 

{\bf One arm events:} 
Let $S$ be a connected  measurable subset of $\bR^2$. For $A,B\subset S$ with $A,B$ connected and $A\cap B=\phi$, 
\begin{align*}
A  \mathrel{\mathop{\llra}^{S}} B:= &\mbox{ the event that there exists a path in the graph $G_\lam^e$ from some point in } A \mbox{ to some point in } B   \nn\\
&\; \mbox{ entirely confined in } S.\nn
\end{align*}
Similarly for connected subsets $C\subset D \subset \bR^2$
\begin{align*}
C \llra \partial D := & \mbox{ the event that there exists a path in the graph $G_\lam^e$ from some point in } C \mbox{ to some point in } D^c.
%&\; \mbox{ entirely confined in }  Q.
\end{align*}
\subsection{Main Results. }
Our first result is on the existence of a phase transition in the three models described earlier. Penrose \citep{Penrose1991} showed that there is a non-trivial phase transition in the RCM, that is, $0<\lam_c<\infty$  under the condition $0<\int_{0}^{\infty} rg(r)dr<\infty$. We now prove a similar result for the eRCM albeit under a stronger restriction on $g$. In case of the ieRCM the condition required is stronger than the one for the iRCM as derived in \citep{Deprez2018}. \citep{RRoy1991} showed the existence of a phase transition in the Poisson stick model under the assumption that the stick length distribution has bounded support, a result which we extend to sticks of unbounded lengths.

\begin{thm}
A phase transition occurs in the
\begin{enumerate}[(i)]
\item eRCM $G_\lam^e$ if the connection function $g$ satisfies $0<\int_0^{\infty}r^3 g(r)\;dr <\infty$.
\lab{np_ercm}
\item ieRCM $H^e_\lam$ with the connection function of the form (\ref{eq:connection_function}) if  $\al>4$ and $\al\be>8. $
\lab{np_iercm}
\item graph $PS_\lam$ with half length density $h$ if  $0<\int_{0}^{\infty}\ell^2\,h(\ell)\,d\ell<\infty$.
\lab{np_ps}
\end{enumerate}
\lab{thm:PhaseTransition}
\end{thm}
The condition for the existence of an infinite component for large intensities is obtained by comparing the enhanced model with the usual non-enhanced version. For the other side we shall show that the probability that there exists a self avoiding path of length $n$ converges to zero as $n \to \infty$.

Our next result establishes a RSW lemma which is one of the most useful result in planar percolation models. It states that if the probability of crossing a square is uniformly bounded away from zero then so is the probability of crossing a rectangle along the longer side. We demonstrate its utility by establishing that percolation does not occur at criticality. The RSW lemma was first proved for the independent Bernoulli bond percolation model on the $\bZ^2$ lattice  independently by Russo \citep{Russo1978} and Seymour and Welsh \citep{SW1978}. A similar result about occupied and vacant crossings was proved for Boolean model on $\bR^2$ by Roy\citep{RRoy1990}. In Roy\citep{RRoy1991} an RSW result has been proved for the Poisson stick model with sticks of bounded lengths. The RSW results in this article are analogous to those in \citep{Tassion2016} for the percolation model on Poisson-Voronoi tessellations in $\bR^2$. RSW results, continuity of critical parameter and sharpness of phase transition for the Boolean model with unbounded radius distribution has been studied by Ahlberg et. al. in \citep{Ahlberg2018}. For the continuum percolation model with random ellipses on the plane, percolation and connectivity behavior of the vacant and covered set has been studied by Teixeira et. al. in \citep{Teixeira2017}. 

We prove the RSW lemma under the condition that the connection function in the eRCM is of the form $g(r) = O(r^{-c})$ as $r\to \infty$ and the half length density $h$ satisfies $h(\ell) = O(\ell^{-c})$  as $\ell\to \infty$. This assumption is required in order to derive an estimate on the longest edge/stick length intersecting a large box. By Theorem~\ref{thm:PhaseTransition} a phase transition occurs under the above assumptions in the eRCM provided $c > 4$ and in $PS_{\lam}$ if $c > 3$.
%$g(x,y)= 1-\exp\left(-\f{\eta W_x W_y}{|x-y|^{\al}}\right)$  
%, for the eRCM, the ieRCM and the Poisson stick models respectively.   %The RSW lemma for the eRCM and the Poisson stick models under assumption on the tail decay of the connection function. 
%
The first two assertions in Theorem~\ref{thm:RSW_nonpercolation} are equivalent formulations of the RSW lemma while the third is a derived as a (non-trivial) corollary of the first assertion. The RSW results are derived by adapting the technique developed by Tassion in \citep{Tassion2016}. 
We then use a renormalization technique similar to the one used in \citep{Alexander1996} (see also \citep{Daniels2016}) to show that the parameter set over which percolation occurs is open. We shall prove the results in detail for the eRCM. Much of the proof carries over to the other two models for which we will provide only the necessary details.

\begin{thm}
Suppose the following conditions hold.
\begin{enumerate}[(I)]
\item In the eRCM $G_\lam^e$ the connection function $g$ satisfies $g(r)=O(r^{-c})$ as $r\rar \infty$ with $c> 4$.
\item In the ieRCM $H_\lam^e$ the connection function $g$ is of the form (\ref{eq:connection_function}) with $\min\{\al,\al\be\}> 4$.
\item In the graph $PS_\lam$ with half length density $h$ satisfies $h(\ell)=O(\ell^{-c})$ with $c> 3$.
\end{enumerate}
Then the following conclusions hold for all the three graphs $G_\lam^e$, $H_\lam^e$ and $PS_\lam$.
\begin{enumerate}[(i)]
\item If $\inf\limits_{s \geq 1} C_s(1)>0$ then for any $\rho \geq  1$, $\inf\limits_{s\geq 1} C_s(\rho)>0.$
\lab{rsw_inf}
\item  If $\lim\limits_{s\rar \infty} C_s(1)=1$ then for any $\rho \geq 1$, $\lim\limits_{s \rar \infty} C_s(\rho)=1.$
\lab{rsw_lim}
\item  The percolation function is continuous.
\lab{non_perc}
\remove{
\begin{equation}
\liminf_{s>0} C_s(\rho)>0.
\lab{eq:limit_inf_C_s1}
\end{equation}
\lab{crossing_rectangle_inf}

\item  If $\lim\limits_{s\rar \infty} C_s(1)=1$ then
\begin{equation}
\lim_{s \rar \infty} C_s(\rho)=1.
\lab{eq:limit_C_s1}
\end{equation}
\lab{crossing_rectangle_exact}

\item  The set of parameters $\lam$ for which percolation occurs is an open set, in particular there is no percolation at criticality.
\lab{nonpercolation_in}}
\end{enumerate}
\lab{thm:RSW_nonpercolation}
\end{thm}
\section{Proofs}
In what follows $c_0,c_1,c_2,\cdots$ and $C_1,C_2,\cdots$ will denote constants whose values will change from place to place. $|\cdot|$ will be used to refer to the Euclidean norm, the cardinality of a set as well as the Lebesgue measure. 

\remove{The condition for the existence of an infinite component for large intensities is obtained by comparing the enhanced model with the usual non-enhanced version. For the other side we will bound the component containing the origin by a sub-critical branching process in case of the eRCM that dies out with probability one. For the ieRCM we evaluate the probability of a self avoiding path of length $n$ and then show that the probability that there is such a path converges to zero as $n \to \infty$. In order to show that the percolation function is continuous, we first derive a RSW Lemma which is interesting in its own right.  We do this by adapting the technique developed in Tassion \citep{Tassion2016}. We then use a renormalization technique similar to the one used in Daniels \citep{Daniels2016} to show that the parameter set over which percolation occurs is open. In this case we will prove the results in detail for the eRCM. Much of the proof carries over to the other two models for which we will provide only the necessary details.}

\subsection{Proof of Theorem~\ref{thm:PhaseTransition}~(\ref{np_ercm})}
It is clear from the definition that $G^e_{\lam}$ percolates  if $G_{\lam}$ does. So we have $\lam^e_c \leq \lam_c$. From Theorem 1 in Penrose \citep{Penrose1991} we know that $\lam_c \in (0, \infty)$ iff $\int_0^{\infty}r g(r)\;dr \in (0, \infty)$. Since $g(r) \in [0,1]$, $\int_0^{\infty}r^3g(r)\;dr  \in (0, \infty)$ implies $\int_0^{\infty}r g(r)\;dr  \in (0, \infty)$. It follows from the above observations that $\lam_c^e < \infty$. 

\begin{figure}[h!]
\centering
\includegraphics[width=1.03\linewidth]{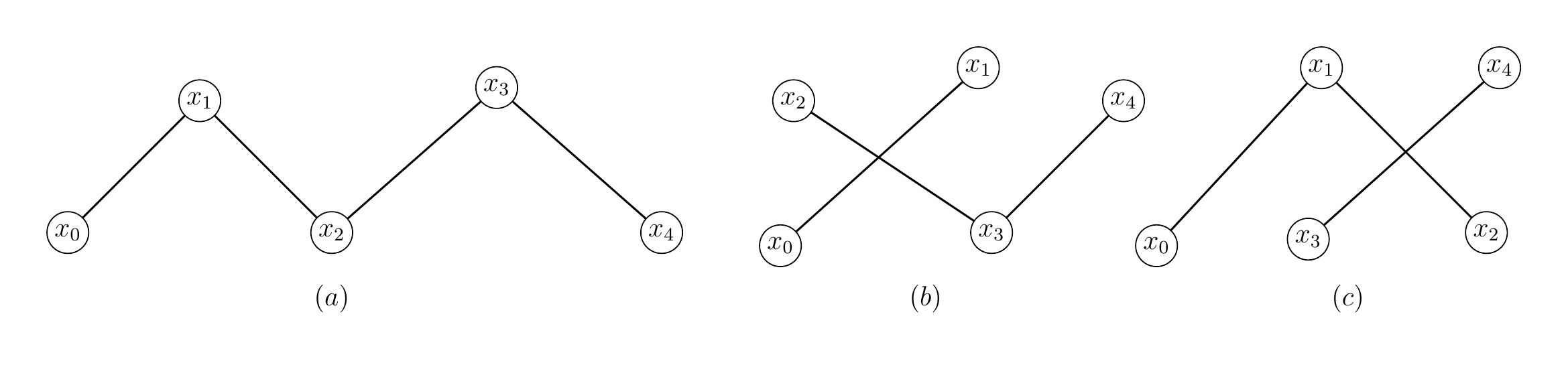}
  \captionof{figure}{Paths from $x_0 \to x_1 \to x_2 \to x_3\to x_4$}
  \label{fig:No_of_Paths}
\end{figure}

We now show that $\lam^e_c> 0$. We shall bound the probability that there is a {\it self-avoiding path} formed using $n$ distinct points of $\cP_{\lam}$ starting from the origin. For any $\textbf{x} = (x_1, x_2, \cdots, x_n)$ of $n$ distinct points of the Poisson point process, we examine whether there is a path starting at $x_0=O$ and uses edges (or parts of it) with end points from $\textbf{x}$ in the coordinate order. We shall denote any such path by $x_0 \to x_1 \to x_2 \to \cdots \to x_n$ even though some of these points may not be part of the path. While the path may have loops, each edge or a part of it is used exactly once while traversing the path. For each ordered sequence $\textbf{x}$ as above, there can be several ways in which a path can occur. See Figure~\ref{fig:No_of_Paths} for self-avoiding paths formed by $(x_1, x_2,x_3,x_4)$. Each such possibility gives rise to a unique {\it block structure} that we describe below. In order to carry out the computation we segregate all paths in disjoint block structures.

We now illustrate this via an example. Take $n=4$, $x_0=O$ and $x_1, x_2, x_3, x_4$ be four distinct points in $\cP_\lam$. Suppose that $O \to x_1 \to x_2 \to x_3 \to x_4$ is a self avoiding path. This can occur in only one way in the RCM (Figure~\ref{fig:No_of_Paths}~(a)) but in the eRCM this can occur in three different ways (see Figure~\ref{fig:No_of_Paths} ). Note that in Figure~\ref{fig:No_of_Paths} (c) we allow the segment $\overline{x_3 x_4}$ to intersect the segment $\overline{x_0 x_1}$.
%Each diagram represents a unique block structure that will be described below. There are $4!$ paths comprising the vertices $\{x_1, x_2, x_3, x_4\}$ with a particular block structure starting at $x_0$.

Let $E(G_\lam)$ denote the edge set of the graph $G_\lam$. Fix $n \in \bN$ and let $\textbf{x} = (x_1,x_2, \ldots, x_n) \in \cP_{\lam,\neq}^n$ be an ordered collection of $n$ distinct points in $\cP_{\lam}$. Define the sub-collection of indices
\begin{eqnarray*}
I(\textbf{x}) & := &\{i \in [n-1] : \{x_{i-1}, x_i\}, \{x_i, x_{i+1}\} \in E(G_\lam)\},\\
J(\textbf{x}) & := & \{ i, i+1: i \in [n-2], \{x_{i-1}, x_i\}, \{x_{i+1}, x_{i+2}\} \in E(G_\lam) ,  \overline{x_{i-1}x_i} \cap \overline{x_{i+1}x_{i+2}} = z_i \not\in \cP_{\lam}\}.
\end{eqnarray*} 
The last condition in the definition of $J(\textbf{x})$ requires that the edges intersect at a point interior to both the edges. Suppose $I(\textbf{x})=\{i_1, i_2,\ldots , i_k\}$ for some $0 \leq k \leq n-1$, labeled in the increasing order. Set $i_0:=0, i_{k+1} := n$ and define the blocks 
\[ B_j(\textbf{x}) := \{i_{j-1} < i < i_j: i \in J(\textbf{x})\} \cup \{i_{j-1}, i_j\}, \qquad 1 \leq j \leq k+1. \]

%For $\textbf{x}:=(x_1,x_2, \ldots, x_n)$ such that 
%We call $P_{\textbf{x}} := \bigcup\limits_{j=1}^{k}\bigcup\limits_{l\in B_j(\textbf{x})} \overline{x_{l}x_{l+1}}$ 
If $x_0 \to x_1 \to x_2 \to \cdots \to x_n$ is a self avoiding path then $\cup_{j=0}^{k+1} B_j(\textbf{x}) = \{0,1,2, \ldots ,n\}$. If $i \in I(\textbf{x})$ then $x_i$ lies on the path whereas if $i \in J(\textbf{x})$ then the path uses a part of an edge one of whose end points is $x_i$. Let $B(\textbf{x}) = (B_1(\textbf{x}), \ldots , B_{k+1}(\textbf{x}))$. With reference to Figure~\ref{fig:No_of_Paths}, the diagram in (a) consists of four blocks $B_i(\textbf{x}) = \{i-1,i\}$, $1 \leq i \leq 4$. The diagram in (b) has two blocks $B_1(\textbf{x}) = \{0,1,2,3\}$, $B_2(\textbf{x}) = \{3,4\}$ while the one in (c) also has two blocks $B_1(\textbf{x}) = \{0,1\}$, $B_2(\textbf{x}) = \{1,2,3,4\}$.

Note that all blocks have even cardinality. For $0 \leq k \leq n-1$, $\cB_k$ be the collection of all block structures $B=(B_1, \ldots , B_{k+1})$ such that $|B_j|$ is even and for some $0 < i_1 <  i_2 < \cdots < i_k < n$, $B_j := \{i: i_{j-1}\leq i \leq i_j\}$ with $i_0:=0, i_{k+1} := n$. Let $B_j^e = \{i_{j-1} + 2r-1: r = 1,2, \ldots , (i_j - i_{j-1} +1)/2\}$ be the set of indices with an even ordering in $B_j$. 

\begin{figure}[h!]
\centering
\includegraphics[width=0.6\linewidth]{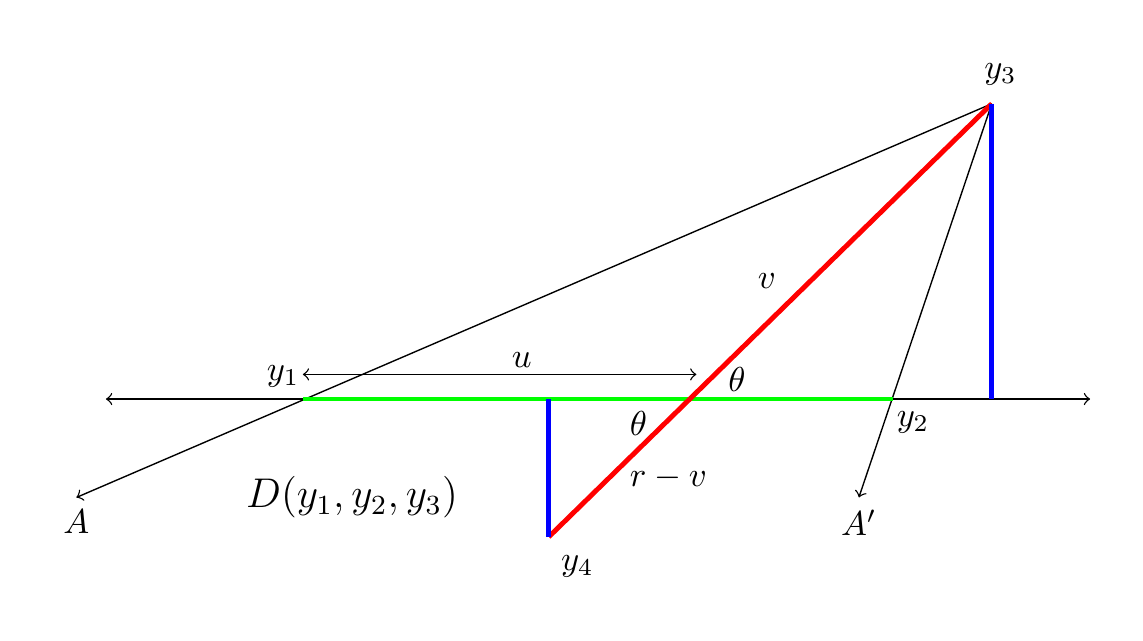}
  \captionof{figure}{$D(y_1, y_2, y_3)$ the unbounded region $A y_1 y_2 A'$}
  \label{fig:region_D}
\end{figure}

By definition of percolation probability the event that the origin lies in an infinite connected component implies that for each $n\in \bN$ there is a self-avoiding path of length $n$ starting from the origin in $G^e_\lam$.
\begin{eqnarray}
\th^e(\lam) &\leq &  P^o(\mbox{ there is a self-avoiding  path on } n \mbox{ vertices in } G_{\lam}^e
\mbox{ starting at } x_0)\nn\\
& \leq &  E^o \left[\sum\limits_{\textbf{x} \in \cP_{\lam,\neq}^n} 1\{x_0 \to x_1 \to x_2 \to \cdots \to x_n \mbox{ occurs}\} \right] \nn\\
& = &   \sum\limits_{k=0}^{n-1}  \sum\limits_{B \in \cB_k}  E^o \left[   \sum\limits_{\textbf{x} \in \cP_{\lam,\neq}^n}  1\{B(\textbf{x}) = B\}\right].
\lab{eq:theta_tilde_2}
\end{eqnarray}
Recall that each $B \in \cB_k$ is of the form $(B_1, \ldots , B_{k+1})$. For a block of the form $B_j := \{i: i_{j-1}\leq i \leq i_j\}$ of size larger than two, let 
\[ A_j = \{\textbf{y}=(y_1,  \ldots , y_n): \overline{y_{i_{j-1} + 2r - 2}y_{i_{j-1} + 2r - 1}} \mbox{ intersects } \overline{y_{i_{j-1} + 2r}y_{i_{j-1} + 2r + 1}} \mbox{ for all } r = 1, \ldots , (i_j - i_{j-1} - 1)/2 \}. \]
The intersection above is understood to be at an interior point of the line segments. For a block of size two we set $1_{A_j} \equiv 1$. Conditioning on $\cP_{\lam}$ and then applying the Campbell-Mecke formula we obtain
\begin{eqnarray}
E^o \left[   \sum\limits_{\textbf{x} \in \cP_{\lam,\neq}^n}  1\{B(\textbf{x}) = B\} \right] & = &  E^o \left[   \sum\limits_{\textbf{x} \in \cP_{\lam,\neq}^n} E^o \left[  \prod\limits_{j=1}^{k+1}1\{B_j(\textbf{x}) = B_j\}\bigg\vert \cP_{\lam} \right]\right]\nn\\
 & = &  E^o \left[   \sum\limits_{\textbf{x} \in \cP_{\lam,\neq}^n} \prod\limits_{j=1}^{k+1} \prod\limits_{l\in B_j^e}g( x_{l-1}, x_l) 1_{A_j}(\textbf{x})  \right] \nn\\
 & = &  \lam^n \int \int \ldots \int  \prod\limits_{j=1}^{k+1} \prod\limits_{l\in B_j^e}g( x_{l-1}, x_l) 1_{A_j}(\textbf{x}) \prod\limits_{l=1}^{n} dx_l. 
\lab{eq:theta_tilde_3}
\end{eqnarray}
%
%Applying Campbell-Mecke formula for the expectation inside the sum on the right hand side of (\ref{eq:theta_tilde_2}) equals
%
%\begin{equation}
% E^o \left[   \sum\limits_{\textbf{x} \in \cP_{\lam,\neq}^n} E^o \left[ \prod\limits_{j=1}^{k+1} \prod\limits_{l\in B_j^e}g( x_{l-1}, x_l) \bigg\vert \cP_{\lam} \right]\right]  
  %
%= \lam^n \int \int \ldots \int  \prod\limits_{j=1}^{k+1} \prod\limits_{l\in B_j^e}g( x_{l-1}, x_l)  \prod\limits_{l=1}^{n} dx_l. 
% =  \lam^n \int \int \ldots \int    \prod\limits_{\stackrel{1 \leq j \leq k+1}{j \mbox{ odd}}} \prod\limits_{l\in B_j^e} g( x_{l-1}, x_l)\prod\limits_{\stackrel{1 \leq j \leq k+1}{j \mbox{ even}}} \prod\limits_{l\in B_j^e} g( x_{l-1}, x_l) \prod\limits_{l=1}^{n} dx_l. 
%
%\lab{eq:theta_tilde_3}
%\end{equation}
%
We now evaluate the contribution to (\ref{eq:theta_tilde_3}) from blocks of various sizes. By contribution from a block we refer to the outcome of evaluating the integrals in (\ref{eq:theta_tilde_3}) with respect to all variables with index in that block except for the one corresponding to the first index. We shall integrate the variables in the descending order starting with those in the block $B_{k+1}$. Blocks of size four and higher yield a nice formula for the upper bound. To see this one needs to compute the bound for a block of size six. We start with the simplest block of size two. The contribution from a block of size two will be of the form
%the form $\{y_1,y_2\}$ is
%
\begin{equation}
\int_{\bR^2}  g( |x_1-x_2|)\, dx_2 = \int_{\bR^2}  g(|x|)\, dx= 2\pi \int_{0}^{\infty} r\, g( r)\, dr.
\lab{eq:theta_tilde_4}
\end{equation} 
Next we compute the contribution from a block of size four. For $a=(a_1,a_2), b=(b_1, b_2)$, let $H^+(a, b) :=\left\{ (c_1,c_2) \in \bR^2 : \f{c_2-a_2}{c_1-a_1} \geq \f{a_2-b_2}{a_1-b_1} \right\}$. For $c\in H^+ (a, b)$, let $D(a, b, c) := \big\{d\in \bR^{2}: \overline{cd} \mbox{ intersects } \overline{ab}\big\}$. The contribution to (\ref{eq:theta_tilde_3}) from a block of size four equals
\begin{equation}
   \int_{\bR^2} g( |y_1-y_2|)\,dy_2\, \int_{H^+(y_1,y_2)}   \int_{D(y_1,y_2, y_3)}g( |y_4-y_3|)\,dy_4 \,dy_3.
\lab{eq:theta_tilde_5}
\end{equation}
The region $D(y_1, y_2, y_3)$ is the region enclosed by the rays $\overrightarrow{y_1A}$, $\overrightarrow{y_2A'}$ and the segment $\overline{y_1y_2}$ (see  Figure~\ref{fig:region_D}). 
%Let $T_x(z)=z+x$ and $R_{\th}$ be the rotation matrix for an angle of rotation $\th$. For fixed $y_1, y_2$, 
Suppose $|y_1-y_2|= \ell$. By a translation and rotation (so that $y_1$ is translated to the origin and $y_2$ to the point $(\ell, 0)$), we can write
\begin{equation}
\int_{H^+(y_1,y_2)}   \int_{D(y_1,y_2, y_3)}g( |y_4-y_3|)\,dy_4 \,dy_3 = \int_{\bR \times \{a_2>0\}} \int_{D(O, (\ell, 0), a)} g(|a-b|)\; db\; da,
\lab{eq:theta_tilde_5a}
\end{equation}
where $a = (a_1,a_2)$ and $b = (b_1, b_2)$. Changing the variables $a, b$ to $u, v, r, \theta$ according to  $a=(u+v \cos\theta, v \sin\theta), b=(u-(r-v)\cos\theta, -(r-v)\sin\theta)$ and noting that the determinant of the Jacobian satisfies $|J^{-1}|=r \sin\theta$ we can rewrite the right hand side of (\ref{eq:theta_tilde_5a}) as
\begin{equation}
% \int_{\bR \times \{a_2>0\}} \int_{D(O, (\ell, 0), a)} g(|a-b|)\; db\; da =   
\int_{0}^{\ell} \int_{0}^{\pi} \int_{0}^{\infty} \int_{0}^{r} g(r)\, r\, \sin\theta \;dv\;dr\;d\theta \;du =  2  \ell \int_{0}^{\infty}r^{2} g(r) \;dr.  
\lab{eq:theta_tilde_5b}
\end{equation}
By changing to polar coordinates and using (\ref{eq:theta_tilde_5b}) the expression in
%(\ref{eq:theta_tilde_5a}) and then using the resulting quantity in 
(\ref{eq:theta_tilde_5}) equals 
\begin{equation}
%\int_{\bR^2} g( |y_1-y_2|)\,dy_2\, \int_{H^+(y_1,y_2)}   \int_{D(y_1,y_2, y_3)}g( |y_4-y_3|)\,dy_4 \,dy_3 
%
\int_{0}^{\infty}\int_{0}^{2\pi}  \ell \, g(\ell)  \left( 2\ell \int_{0}^{\infty}r^{2} \, g(r) \, dr\right) \, d\th \, d\ell = 4\pi  \left( \int_{0}^{\infty}r^{2} \, g(r) \, dr\right)^2.  
\lab{eq:theta_tilde_5c}
\end{equation}
%
%Consider the block $\{y_1,y_2,\ldots, y_6\}$ where  $\overline{y_3 y_4}$ intersects $\overline{y_1 y_2}$ and  $\overline{y_5 y_6}$ intersects $\overline{y_3 y_4}$. 
The contribution from a block of size six can be found by applying the above procedure twice. Using (\ref{eq:theta_tilde_5a}), (\ref{eq:theta_tilde_5b}) we get
\begin{eqnarray}
\lefteqn{  \int_{\bR^2}  g( |y_1-y_2|)\, dy_2 \int_{H^+(y_1,y_2)}   \int_{D(y_1,y_2, y_3)}g( |y_4-y_3|)\,dy_4 \,dy_3  \int_{H^+(y_3,y_4)}   \int_{D(y_3,y_4, y_5)}g( |y_6-y_5|)\,dy_6 \,dy_5}\nn\\
& \leq &  \int_{\bR^2}  g( |y_1-y_2|)\, dy_2 \int_{H^+(y_1,y_2)}   \int_{D(y_1,y_2, y_3)}g( |y_4-y_3|)\,dy_4 \,dy_3 \; 2 |y_4-y_3|  \left(\int_{0}^{\infty} r^2\, g( r)\, dr\right)\nn\\
 & = &   2     \int_{\bR^2} g( |y_1-y_2|)\, dy_2 \; 2 |y_2-y_1|  \left(\int_{0}^{\infty} r^3\, g(r)\, dr\right) \left(\int_{0}^{\infty} r^2\, g(r)\, dr\right)\nn\\
 & = &  8\pi       \left(\int_{0}^{\infty} r^2\, g( r)\, dr\right)  \left(\int_{0}^{\infty} r^3\, g(r)\, dr\right) \left(\int_{0}^{\infty} r^2\, g( r)\, dr\right).
\lab{eq:theta_tilde_6}
\end{eqnarray} 
By iterating the above procedure the contribution from the block of size $2m+2$ on the vertices $\{y_1,y_2,\ldots, y_{2m+2}\}$ where the edge $\overline{y_i,y_{i+1}}$ intersects $\overline{y_{i+2},y_{i+3}}$ for all $i = 1, 2, \ldots , 2m-1$ can be shown to be
\begin{equation}
 2^{m+1}\pi \left(\int_{0}^{\infty} r^2\,  g( r)\, dr\right)  \left(\int_{0}^{\infty} r^3\,  g( r)\, dr\right)^{m-1} \left(\int_{0}^{\infty} r^2\,  g( r)\, dr\right).
 \lab{eq:theta_tilde_6a}
 \end{equation}
 For a path with blocks $B_1, B_2, \ldots , B_{k+1}$ let 
 \begin{equation}
 k_p:=|\{ j : |B_j|=2p \}| \mbox{ and } \bar{k}= \sum\limits_{\stackrel{1 \leq j \leq n}{ |B_j| \geq 6 }} \left( \f{|B_j|}{2}-2 \right).
\lab{eq:Kp}
\end{equation}
  Since each vertex can be in at most two adjacent blocks we have that 
\begin{equation}
\bar{k} \leq \sum\limits_{j=1}^{k+1}\f{|B_j|}{2} \leq n.
\lab{eq:k4_eRCM}
\end{equation}
Substituting the contributions from each block $B_1,B_2, \ldots, B_{k+1}$, that is, the expression in (\ref{eq:theta_tilde_4}), (\ref{eq:theta_tilde_5c})--(\ref{eq:theta_tilde_6a}) in the expression on the right in (\ref{eq:theta_tilde_3}) yields
\begin{equation}
E^o \left[ \sum_{\textbf{x} \in \cP_{\lam,\neq}^n}  1\{B(\textbf{x}) = B\} \right] = \lam^n \prod_{m=0}^{\lfloor \f{n}{2}\rfloor} (2^{m+1}\pi)^{k_m}  \left(\int_{0}^{\infty} r g(r) dr\right)^{k_1} \left(\int_{0}^{\infty} r^2  g(r) dr\right)^{2(k-k_1)} \left(\int_{0}^{\infty} r^3  g(r) dr \right)^{\bar{k}}.
\lab{eq:theta_tilde_7}
\end{equation}
Using (\ref{eq:k4_eRCM}) in  (\ref{eq:theta_tilde_7})  and then using (\ref{eq:theta_tilde_3}), (\ref{eq:theta_tilde_7})  in (\ref{eq:theta_tilde_2})  we obtain
\begin{equation}
\th^e(\lam) \leq \lam^n \sum\limits_{k=0}^{n-1} \sum\limits_{B \in \cB_k} \left( \prod\limits_{m=0}^{\lfloor \f{n}{2}\rfloor} (2^{m+1}\pi)^{k_m} \right)  C_1^{4n} ,
\lab{eq:theta_tilde_8}
\end{equation}
where $C_1 = \max\left\{ \int_{0}^{\infty} r^j\,  g(r)\, dr, j=1,2,3\right\}$. Since $\sum\limits_{m=1}^{\lfloor \f{n}{2}\rfloor}(m+1)k_m \leq \sum\limits_{m=1}^{\lfloor \frac{n}{2}\rfloor}2m\,k_m \leq 2n$ and $|\cB_k| = {n-1 \choose k}$, there exists a constant $C$ such that
\[ \th^e(\lam) \leq  (C \lam)^n \to 0 \]
as $n \to \infty$ for all $\lam \in (0, \lam_0)$, for some $\lam_0 > 0$ provided $\max\left\{ \int_{0}^{\infty} r^j\,  g(r)\, dr, j=1,2,3\right\} <\infty$. \qed
\subsection{Proof of Theorem~\ref{thm:RSW_nonpercolation}~(\ref{rsw_inf}) for eRCM}
\lab{subsection:RSW_for_eRCM}

\begin{figure}
\centering
\includegraphics[width=0.6\linewidth]{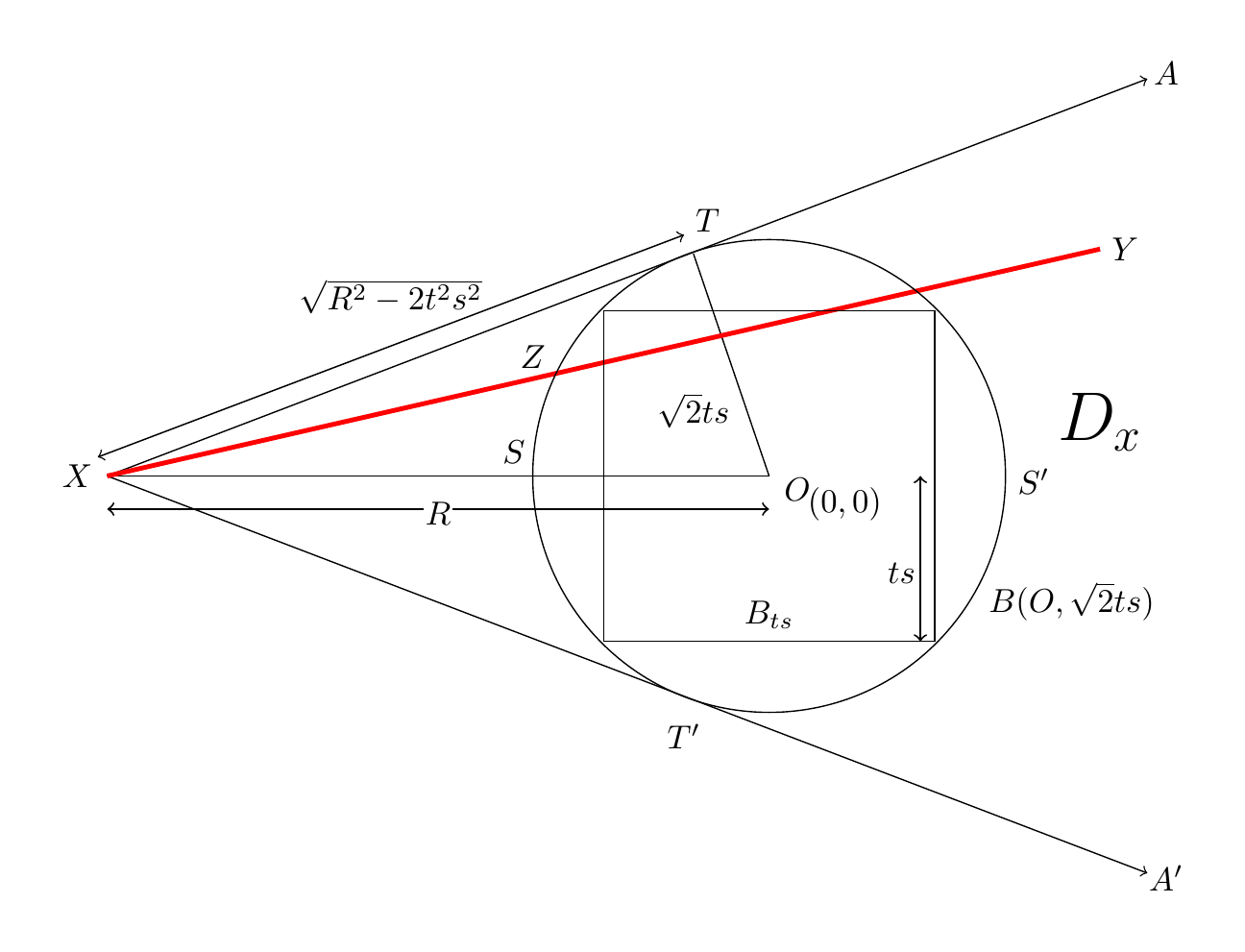}
  \captionof{figure}{$D_x$ the unbounded region $ATS'T'A'$}
  \label{fig:both_out}
\end{figure}

Consider the graph $G_\lam^e$ with connection function $g(r) = O(r^{-c})$ as $r \to \infty$ where $c > 4$ is arbitrary. Later we shall borrow some key results from Tassion \citep{Tassion2016} and hence will adopt many of the notations from that paper. A key ingredient in the proof is the following result on the length of the longest edge in $G_\lam$ which allows us to localize the analysis. 
\begin{prop}
For any $s > 0$ let $M_{s}$ be the length of the longest edge in $G_{\lam}$ intersecting the box $B_{s} = [-s,s]^2$. Suppose that the connection function $g$ satisfies $g(r)=O(r^{-c})$ as $r\rar \infty$. Then for any $c > 4$, $t>0$ and $ \tau>\f{2}{c-2}$ we have $P\left(M_{ts}>s^{\tau}\right)\rar 0$ as $ s\rar \infty$.
\lab{prop:largest_edge_length_in_box}
\end{prop}
{\bf Proof of Proposition~\ref{prop:largest_edge_length_in_box}.}   Fix $c > 4$, $t>0$. Let $B(O,s):=\{x\in \bR^2:|x|\leq s\}$ be the ball of radius $s$ centered at the origin. Recall that for any two points $x,y \in \bR^2$, $\overline{xy}$ denotes the line segment joining $x$ and $y$. Define the events $D_s(l)=\{M_s>l \}$,
$$O_{s}(\tau)=\left\{ X \in \cP_{\lam}: \mbox{ there is an edge of length longer than } s^{\tau} \mbox{ incident on } X \mbox{ in } G_{\lam} \right\}$$ 
and 
$$\bar{O}_{t,s}(\tau)=\left\{ (X, Y) \in \cP_{\lam}^2: \mbox{ there is an edge in } G_{\lam} \mbox{ joining } X,Y,  |\overline{X Y} | \geq s^{\tau}, \overline{X Y} \mbox{ intersects } B(O, \sqrt{2}s)\right\}.$$ 
\begin{eqnarray}
P\left(D_{ts}(s^{\tau})\right)&\leq & E\left[\sum_{X,Y\in \cP_{\lam}}1\left\{\overline{X Y}\mbox{ intersects } B_{ts}\right\}1\left\{|X- Y|\geq s^\tau\right\}\right]\lab{eq:1st_inequality}\nn\\
%& = &  E \left[\sum_{X,Y\in \cP_{\lam}}1_{\{\overline{X Y} \mbox{ intersects } B_{ts}\}}1_{\{\mbox{at least one of } X,Y \mbox{ is in } B_{ts}\}}1_{\{|X- Y|\geq s^\tau\}}\right] \nn\\
%
%&+& E \left[\sum_{X,Y\in \cP_{\lam}}1_{\{\overline{X Y} \mbox{ intersects } B_{ts}\}}1_{\{\mbox{neither of } X,Y \mbox{ is in } B_{ts}\}}1_{\{|X- Y|\geq s^\tau\}}\right] \nn\\
& \leq &  E\left[\sum_{X\in \cP_{\lam}\cap B(O,\sqrt{2}ts)}1\left\{X\in O_{s}(\tau)\right\}\right] + 
E \left[\sum_{X,Y \in \cP_{\lam}\cap B(O,\sqrt{2}ts)^c }1\left\{(X,Y) \in \bar{O}_{t,s}\right\}\right].
\lab{eq:2nd_inequality}
\end{eqnarray}
The Campbell-Mecke formula applied to the first term on the right hand side of the last inequality in (\ref{eq:2nd_inequality}) yields
\begin{eqnarray}
E\left[\sum_{X\in \cP_{\lam}\cap B(O,\sqrt{2}ts)}1\left\{X\in O_{s}(\tau)\right\}\right]
&=& C \lam (ts)^{2}P^o\left(O\in O_{s}(\tau)\right)\nn\\
&=& C \lam (ts)^{2}\left(1-P^o\left(\mbox{none of the edges incident on O is of length} \geq s^{\tau}\right)\right)\nn\\
& = & C \lam (ts)^2\left(1- \exp\left(-\lam \int_{B(O, s^{\tau})^c}g(|x|)dx\right)\right)\nn\\
 & \leq & C (\lam ts)^2  \int_{B(O, s^{\tau})^c}g(|x|)dx =C_1 s^2 \int_{s^\tau}^{\infty}rg(r)\,dr \leq C_2\, s^{2-\tau(c-2)},
\lab{eq:one_in}
\end{eqnarray}
where we have used the fact that the points of $\cP_{\lam}$ from which there is incident on $O$ an edge that is of length longer than $s^{\tau}$ is a Poisson point process of intensity $\lam g(|x|) 1\{x \in B(O, s^{\tau})^c\}$, the inequality $1-e^{-y}\leq y$ and the assumption on $g$. Similarly we can bound the second term on the right hand side in the last inequality in (\ref{eq:2nd_inequality}) as follows (see Figure~\ref{fig:both_out}).
\begin{equation}
 E \left[\sum_{X,Y \in \cP_{\lam}\cap B(O,\sqrt{2}ts)^c }1\left\{(X,Y) \in \bar{O}_{t,s}\right\}\right] = \lam^2 \int_{ B(O,\sqrt{2}ts)^c}\int_{ D_x \cap B(x, s^{\tau})^c}g(|x-y|)\,dy\,dx,
\lab{eq:both_out_1}
\end{equation}
where $D_x$ is the unbounded region $ATS'T'A'$ as shown in Figure~\ref{fig:both_out}. Changing to polar coordinates and using the obvious bounds for the range of the $y$-variable we can bound the expression on the right in (\ref{eq:both_out_1}) by
\begin{eqnarray}
\lefteqn{ C_3 \int_{\sqrt{2}ts}^{\infty}\int_{ 0}^{2\pi}R\,d\phi\,dR\int_{s^\tau \vee \sqrt{R^2-2t^2s^2}} ^{\infty}rg(r)\,dr}\nn\\
&\leq& C_4 \int_{\sqrt{2}ts}^{\infty}\left(s^\tau \vee\sqrt{R^2-2t^2s^2}\right)^{2-c}R\,dR\nn\\
&=& C_4 \int_{ ts}^{\sqrt{2t^2s^2 + s^{2\tau}}}s^{\tau(2-c)}R\,dR + C_4 \int_{\sqrt{2t^2s^2+s^{2\tau}}}^{\infty} \left(R^2-2t^2s^2 \right)^{\f{2-c}{2}}R\,dR\nn\\
&=& C_5 s^{-\tau(c-4)} + C_6 \int_{s^\tau}^{\infty} u^{3-c}\,du = C_7 s^{-\tau(c-4)},
\lab{eq:both_out_2}
\end{eqnarray}
where we have used the assumption that $g(r)\leq C r^{-c}$. 
Substituting from (\ref{eq:one_in}) and (\ref{eq:both_out_2})  in (\ref{eq:1st_inequality}) we obtain
\begin{eqnarray} 
P\left(D_{ts}(s^\tau)\right)&\leq&  C_2 \, s^{2-\tau(c-2)} + C_7 s^{-\tau(c-4)}.\nn
\end{eqnarray}
Hence $P\left(D_{ts}(s^{\tau})\right) \rar 0$ as $s\rar \infty$, since $\tau>\f{2}{c-2}$ and $c>4$.\qed

The following corollary gives us the precise form in which we will be using Proposition~\ref{prop:largest_edge_length_in_box}.
\begin{cor} 
For the graph $G_\lam$ with the connection function $g$ satisfying $g(r)=O(r^{-c})$ as $r\rar\infty$, let $L_{ts}(\tau)$ be the event that there exists an edge of length longer than $s^\tau$ intersecting the annulus $A_{2ts,4ts}$. Then for any $c > 4$, $t>0$ and $\tau > \f{2}{c-2}$ we have $P\left(L_{ts}(\tau)\right)\rar 0$ as $ s\rar \infty$.
\lab{cor:edge_annulus}
\end{cor}
By assumption, we have for some $c_0 > 0$,
\begin{equation}
C_s(1)\geq c_0\mbox{ for all } s\geq 1.
\lab{eq:assumption_1}
\end{equation} 
Proposition~\ref{prop:general_s_2} below is a restatement of the first assertion in Theorem \ref{thm:RSW_nonpercolation} for the case $\rho = 2$. We shall first use this proposition to extend the result for general $\rho$ and then follow it up with the proof of the proposition. 
\begin{prop}
%In $G_\lam^e$ with the connection function $g$ satisfying $g(r)=O(r^{-c})$ as $r\rar \infty$ with $ c>4$, 
Suppose (\ref{eq:assumption_1}) holds for the graph $G_{\lam}^e$ with the connection function $g$ satisfying $g(r)=O(r^{-c})$ as $r\rar \infty$ for some $ c>4$. Then $\inf\limits_{s \geq 1} C_s(2) > 0.$
\lab{prop:general_s_2}
\end{prop}
\begin{figure}[h!]
  \centering
\hspace{-0.6in}
\includegraphics[width=0.95\linewidth]{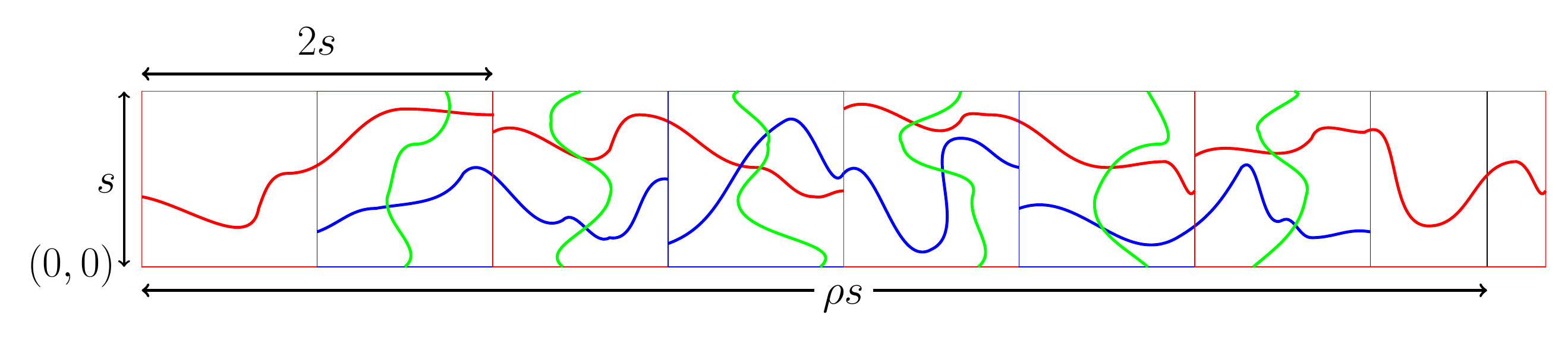}
  \captionof{figure}{A realization of the event $F_s(\rho)$}
  \label{fig:overlapping_paths_1}
\end{figure}
Let $s \geq 1$. Assuming that Proposition~\ref{prop:general_s_2} holds, it suffices to prove the result for $\rho > 2$. We need to build a left to right crossing in  $[0,\rho s]\times[0,s]$. Observe that 
\begin{equation}
[0,\rho s]\times[0,s]\subset \bigcup\limits_{j=0}^{n_\rho}\left(js,0)+[0,2s]\times[0,s]\right),
\lab{eqn:bound_LR_Crossing}
\end{equation}
where $n_\rho\leq \lfloor \rho \rfloor$. Let $F_s(\rho)$ be the event that there exists left to right crossing in $(js,0)+[0,2s]\times[0,s]$ for all $j=0,1,2,\cdots, n_\rho$ and top to down crossing in $(js,0)+[0,s]\times[0,s]$ for all $j=1,2,\cdots, n_\rho$. From (\ref{eqn:bound_LR_Crossing}) we have $F_s(\rho)\subset LR_s(\rho)$ (see Figure~\ref{fig:overlapping_paths_1}). Using this inclusion and applying the FKG inequality we obtain
\[ C_s(\rho)\geq P[LR_s(2)]^{n_\rho+1}P[TD_s(1)]^{n_\rho}. \] 
The assertion in Theorem~\ref{thm:RSW_nonpercolation}~(\ref{rsw_inf}) now follows from (\ref{eq:assumption_1}) and Proposition~\ref{prop:general_s_2}.

\begin{figure}[h!]
\begin{minipage}{.55\textwidth}
  \centering
  \includegraphics[width=.75\linewidth]{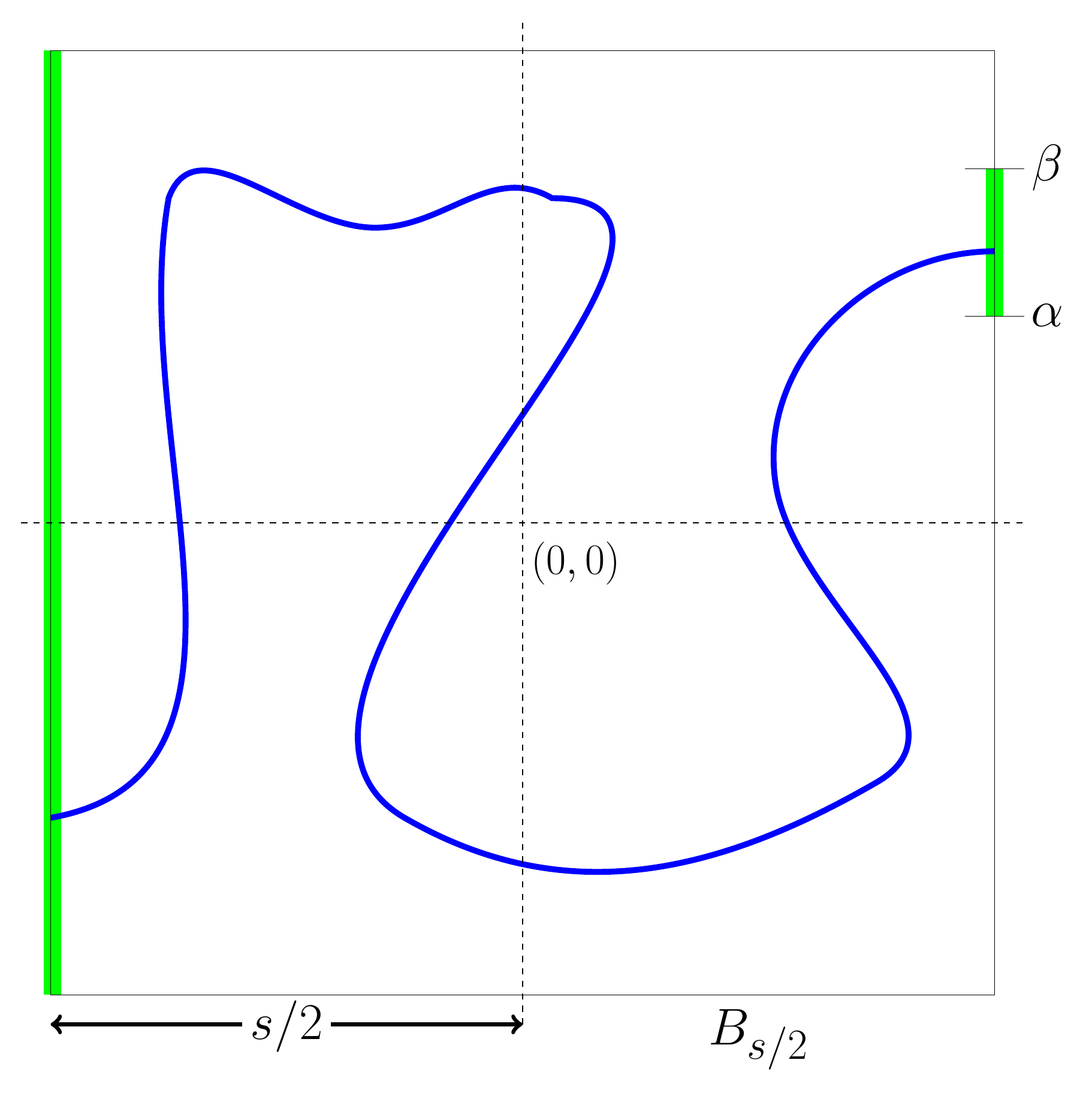}
  \captionof{figure}{Event $\cH_s(\al,\beta)$}
  \label{fig:one_interval}
\end{minipage}%
\begin{minipage}{.55\textwidth}
  \centering
  \includegraphics[width=.75\linewidth]{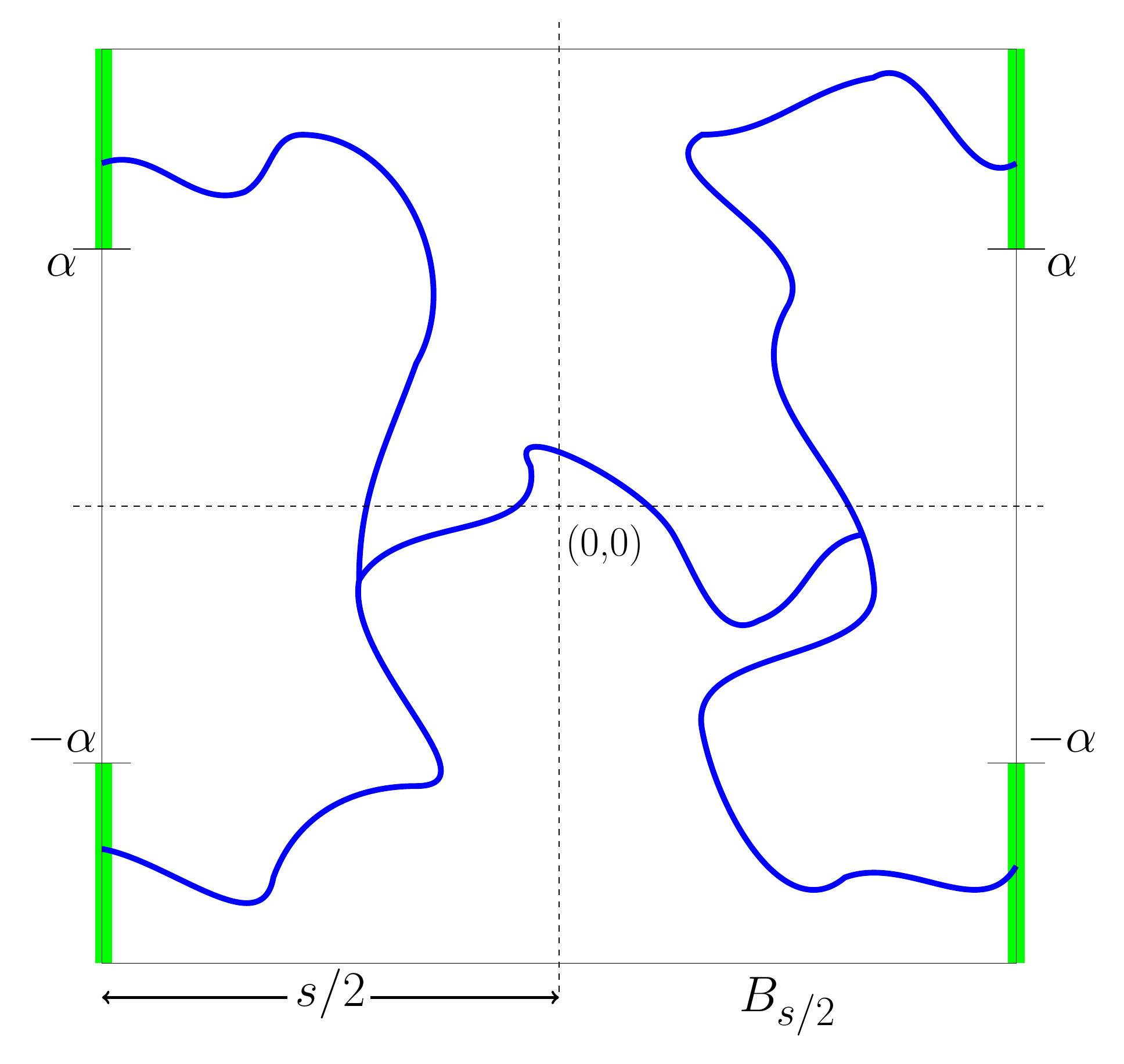}
  \captionof{figure}{Event $\chi_s(\al)$}
  \label{fig:four_interval}
\end{minipage}
\end{figure}

It remains to prove Proposition~\ref{prop:general_s_2}. The proof is derived from the next proposition that follows from a geometric construction. Recall that $\cA_s$ is the event that there exists a circuit in the annulus $A_{s,2s}$.
\begin{prop}
Suppose the conditions given in Proposition~\ref{prop:general_s_2} hold. Then there exists constants $c_2 > 0, C>4$ and an increasing sequence of scales $\{s_n\}_{n\geq 1}$ satisfying $4s_n\leq s_{n+1} \leq Cs_n$ such that  $P[\cA_{s_n}]\geq c_2$ for all $n\geq 1$.
\lab{prop:sequence_of_good_scales}
\end{prop}
{\bf Proof of Proposition~\ref{prop:sequence_of_good_scales}.}
Fix $s\geq 1$.  For $\al, \beta \in [-s/2,s/2],\,  \al<\be$, let $\cH_s(\al,\beta)$ be the  event (see Figure~\ref{fig:one_interval}) that there exists a path in the box $B_{s/2}$ from left to $\{s/2\}\times [\al,\beta]$. For $\al \in[0,s/2]$, define $\chi_s(\al)$ be the event that there exists a path from $\{-s/2\}\times [-s/2,-\al]$ to  $\{s/2\}\times [\al,s/2]$ and there exists a path from $\{-s/2\}\times [\al, s/2]$ to $\{s/2\}\times [-\al,-s/2]$ in $B_{s/2}$ (see Figure~\ref{fig:four_interval}). 
%We call the event $\chi_s(\al)$ by \textit{fork}. 
%
Let  $\cH'_s(\al,\be)$ be the event that there exists a path in the box $B_{s/2}$ from right to $\{-s/2\}\times [\al, \be]$ and $TD'_{s}(1)$ be the event that there exists a top to down path in the box $B_s$. By rotation and reflection invariance the events $\cH_s(-\al, -s/2) ,  \cH'_s(\al,-s/2),  \cH'_s(-\al,-s/2)$ are symmetric versions of $\cH_s(\al,s/2)$. Observe that 
\[ \cH_s(\al,s/2) \cap \cH_s(-\al,-s/2) \cap \cH'_s(\al,s/2) \cap \cH'_s(-\al,-s/2) \cap TD'_{s/2}(1) \subset \chi_s(\al).\]

Given the assumption (\ref{eq:assumption_1}) that the probability of crossing boxes $C_s(1)$ is uniformly bounded away from zero, the following two Lemmas from Tassion\citep{Tassion2016} provide lower bounds for probabilities of certain paths that will allow us to glue them together to construct paths with desired properties. These Lemmas are true for any planar percolation model which in conjunction with a result such as Corollary~\ref{cor:edge_annulus} allows us to derive the RSW result.
\begin{lem} 
If for some $c_0>0$, $\inf_{s \geq 1}C_s(1)\geq c_0$, then for all $s\geq 1$ there exists $ \al_s \in [0,s/4]$  and $c_1>0$ such that,
\begin{enumerate}[(i)]
\item for all $0\leq \al \leq \al_s$
\begin{equation}
P\left(\chi_s(\al)\right)\geq c_1,
\lab{eq:prob_chi}
\end{equation} \label{part1_chi}
\item if $\al_s<s/4$, then for all $\al_s \leq \al \leq s/2 $
\begin{equation}
 P\left(\cH_s(0,\al)\right)-P\left(\cH_s(\al,s/2)\right)\geq c_{0}/4.
\lab{eq:prob_phi}
\end{equation}\label{part2_phi}
\end{enumerate}
\lab{lem:lowerbound_on_chi&phi}
\end{lem}
\begin{lem}
Let $\al_s$ be as in Lemma~\ref{lem:lowerbound_on_chi&phi}. Then the following two statements are true.
\begin{enumerate}[(i)]
\item There exists $c_2>0$ such that whenever $\al_s\leq 2\al_{\f{2}{3}s}$ for some $s\geq 2$, then 
\begin{equation}
P(\cA_s)\geq c_{2}.
\lab{eq:good_scale_annulus}
\end{equation} \label{part1_goodscales}
\item For any $s\geq 1$ if $ P(\cA_s)\geq c_2$ and $\al_t\leq s$ for some $t\geq 4s$, then there exists $c_3>0$ such that
\begin{equation}
P(\cA_t)\geq c_3. 
\lab{eq:advancement_s_to_t}
\end{equation} \label{part2_s_to_t}
\end{enumerate}
\lab{lem:advancement_s_to_t}
\end{lem}
%
%Assuming the Lemma~\ref{lem:lowerbound_on_chi&phi} and \ref{lem:advancement_s_to_t} we proceed to prove the Proposition \ref{prop:sequence_of_good_scales}. The construction of the sequence $\{s_n\}_{n\geq 1}$ is done  iteratively in the same way as in  Lemma 3.3 of [Tassion'16]. 
Let $c_3$ be as in Lemma~\ref{lem:advancement_s_to_t}~(\ref{part2_s_to_t}) and $c_0>0$ be as in (\ref{eq:assumption_1}). Since $c_0, c_3 \in (0,1)$ and $c > 4$, we can and do choose $C_1 > 16$ such that, 
\begin{equation}
\left(1-\f{c_3}{2}\right)^{\lfloor\log_5 \f{C_1}{2}\rfloor-1}< \f{c_0}{4}
\lab{eq:condition_on_constants}
\end{equation}
and $\tau \in \left(\f{2}{c-2}, 1\right)$.
%Recall the definition of $L_{ts}(\tau)$, the event that there exists an edge of length larger than $s^\tau$ with at least one end vertex inside the annulus $A_{2ts,4ts}$, for some fixed $t>0$. \red{The definition of $L_{ts}(\tau)$ is a bit different in the Corollary.} 
%Recall that the connection function $g$ satisfies $g(r)=O(r^{-c})$ as $r\rar\infty$ for some $c > 4$. Fix $\tau \in (\f{2}{c-2}, 1)$, the interval being nonempty since $c > 4$. 
%
By Corollary~\ref{cor:edge_annulus} there exists $s_0 \geq 1$ such that 
%for all $ i=2,3,\ldots, \lfloor\log_5 \f{C_1}{2}\rfloor$ and for all $s\geq s_0$
%
\begin{equation}
 P\left(L_{\f{5^i}{4}s}(\tau) \right)\leq \f{c_3}{2} \qquad \mbox{for all }  \; i=2,3,\ldots, \lfloor\log_5 \f{C_1}{2}\rfloor \mbox{ and } s \geq s_0.
\lab{eq:measurable_L}
\end{equation}
Let $\al_s$ be as in Lemma~\ref{lem:lowerbound_on_chi&phi}. Since $\al_s < s$ there must exist a $s_1 >  s_0$ such that
%The sublinear growth of $\al_s$ implies the existence of $s_1\geq s_0$ such that
% 
\begin{equation}
\al_{s_1}\leq 2\al_{\f{2}{3}s_1}.
\lab{eq:goodscales_s1}
\end{equation}
By  (\ref{eq:goodscales_s1}) and Lemma~\ref{lem:advancement_s_to_t}(\ref{part1_goodscales}) we have 
\begin{equation}
P(\cA_{s_1})\geq c_2.
\lab{eq:annulus_s_1}
\end{equation}
Having found $s_1$ the next task is to find $s_2$. This is done using the two steps described in the following Lemma.
\begin{lem}
Let $C_1$ satisfy (\ref{eq:condition_on_constants}) and $s_0 \geq 1$ be such that (\ref{eq:measurable_L}) holds. If $P(\cA_{s})\geq c_2$ for any $s \geq s_0$, then there exists $s'\in[4s,C_1s]$ such that $\al_{s'}\geq s$. 
Further, there exists a constant $C_1'$ and  $s_2 \in [s', C_1's']$ such that $\al_{s_2}\leq 2\al_{\f{2}{3}s_2}$. %where $C_1'=(C_1/3)^{\log_{\f{4}{3}}(\f{3}{2})}$.
\lab{lem:condition_on_constants}
\end{lem}
We now complete the proof of Proposition~\ref{prop:sequence_of_good_scales}.
By Lemma \ref{lem:condition_on_constants} and  (\ref{eq:annulus_s_1}) there exists $s_1' \in  [4s_1, C_1s_1]$ such that $\al_{s'_1}\geq s_1$.
%
%\begin{equation}
%al_{s'_1}\geq s_1.
%\lab{eq:alpha_s1'}
%\end{equation} 
%
Consequently by the second assertion Lemma~\ref{lem:condition_on_constants} 
%and (\ref{eq:alpha_s1'})  
there exists a $C_1'$ and $s_2 \in [s_1', C_1's_1']$ such that $\al_{s_2}\leq 2\al_{\f{2}{3}s_2}$. 
%where $C_1'=(C_1/3)^{\log_{\f{4}{3}}(\f{3}{2})}$. 
An application of Lemma~\ref{lem:advancement_s_to_t}~(\ref{part1_goodscales}) now yields
\begin{equation}
P(\cA_{s_2})\geq c_2.
\lab{eq:s1_to_s2}
\end{equation}
Observe that $4s_1\leq s_1' \leq s_2\leq C_1's_1'\leq C_1C_1's_1$. Setting $C=C_1C_1'$
%=C_1(C_1/3)^{\log_{\f{4}{3}}(\f{3}{2})}$.  
and iterating this procedure we obtain the desired sequence $\{s_n\}_{n\geq 1}$. This proves the Proposition~\ref{prop:sequence_of_good_scales} except for Lemma~\ref{lem:condition_on_constants}.

{\bf Proof of Lemma~\ref{lem:condition_on_constants}.}
The first part of the proof of Lemma~\ref{lem:condition_on_constants} is similar to Lemma 3.2 in \citep{Tassion2016}.
Suppose that for some $s\geq s_0$ and $c_2>0,\, P(\cA_s)\geq c_2$. Suppose if possible $\al_t<s$ for all   $t \in [4s,C_1s]$. If we take $ t=C_1s $, then this yields $\al_{C_1 s} < s < \f{C_1}{4}s$.  It follows by Lemma~\ref{lem:lowerbound_on_chi&phi}~(\ref{part2_phi}) that
\begin{equation}
P\left(\cH_{C_1s}\left(0,s\right)\right)-P\left( \cH_{C_1s}\left(s,\f{C_1s}{2}\right)\right)\geq \f{c_0}{4}.
\lab{eq:prob_chi_with_constant}
\end{equation}
We will now derive a contradiction to (\ref{eq:prob_chi_with_constant}). Note that $\f{5^is}{2} \in [4s, C_1s]$ for $i=2,3, \ldots , \lfloor \log_5\f{C_1}{2} \rfloor$. Since $P(\cA_s) \geq c_2$, taking $t=\f{5^is}{2}$ we have by Lemma~\ref{lem:advancement_s_to_t}~(\ref{part2_s_to_t})
\begin{equation}
 P\left(\cA_{\f{5^i}{2}s}\right) \geq c_3.
\lab{eq:prob_scaled_annulus}
\end{equation}
Fix $\tau < 1$ be such that (\ref{eq:measurable_L}) holds for all $s\geq s_0$. Combining (\ref{eq:prob_scaled_annulus}) and (\ref{eq:measurable_L}) we can write for  $i=2,3,\cdots, \lfloor\log_5 \f{C_1}{2}\rfloor$,
\begin{equation}
 P\left(\cA_{ \f{5^i}{2}s}\cap L^c_{\f{5^i}{4}s}(\tau)\right) \geq \f{c_3}{2}.
\lab{eq:combined_A_L}
\end{equation}
Consider $\cE_s$ to be the event that there exists a circuit in $A_{s,\f{C_1s}{2}}$. Observe that if $\cA_{ \f{5^i}{2}s}\cap L^c_{\f{5^i}{4}s}(\tau)$ occurs for some $ i=2,3,\cdots, \lfloor\log_5 \f{C_1}{2}\rfloor$,  then $\cE_s$ will also occur. Hence
\begin{equation}
 P(\cE_s)\geq P\left(\bigcup_{i=2}^{\lfloor\log_5\f{C_1}{2}\rfloor}  \left(\cA_{ \f{5^i}{2}s}\cap L^c_{\f{5^i}{4}s}(\tau)\right)\right).
\lab{eq:event_E_0}
\end{equation}
In $G^e_\lam$  for any $A_{2s,4s}-$ measurable event $E$, the event $E\cap L_s(\tau)^c$ is measurable with respect to the Poisson point process $\cP_{\lam}$ restricted to the region $A_{s, 5s}$ for $\tau<1$. Indeed, if $E$ is an event measurable w.r.t. $A_{2s,4s}$ and if  the event  $L_{s}(\tau)^{c}$ occurs then there is no edge that intersects $A_{2s,4s}$ and has at least one end vertex out side $A_{2s-s^\tau,4s+s^\tau}$. Hence $E\cap L_{s}(\tau)^{c}$ depends on the Poisson point process $\cP_{\lam}$ restricted to  $ A_{2s-s^{\tau}, 4s+s^{\tau}}$ and $A_{2s-s^{\tau}, 4s+s^{\tau}} \subset A_{s,5s}$ since $\tau < 1$. It follows that $E\cap L_s(\tau)^c$ depends on the Poisson point process $\cP_{\lam}$ restricted to $A_{s, 5s}$. Consequently the events  $\cA_{ \f{5^i}{2}s}\cap L^c_{\f{5^i}{4}s}(\tau)$ for $ i=2,3,\cdots, \lfloor\log_5 C_1\rfloor$ are independent. Using this fact in (\ref{eq:event_E_0}) and  by substituting from (\ref{eq:combined_A_L}) and (\ref{eq:condition_on_constants}) yields
\begin{eqnarray}
P(\cE_s^c)&\leq& \prod\limits_{i=2}^{\lfloor\log_5 \f{C_1}{2}\rfloor} P\left( \left(\cA_{ \f{5^i}{2}s}\cap L^c_{\f{5^i}{4}s}(\tau)\right)^{c}\right) \nn \\
&<& \left(1-\f{c_3}{2}\right)^{\lfloor\log_5\f{C_1}{2}\rfloor-1}< \f{c_0}{4}.
\lab{eq:event_E_2}
\end{eqnarray}
In the square $B_{C_1s/2}+(-C_1s/2,0)=[-C_1s,0]\times [-C_1s/2,C_1s/2]$ consider the following two events, $LH_s$ is the event that there is a path from left to $\{0\}\times [0,s]$ in $B_{C_1s/2}+(-C_1s/2,0)$ and $UH_s$ is the event that there is path from left to $\{0\}\times [s,C_1s/2]$ in $B_{C_1s/2}+(-C_1s/2,0)$. Observe that when $LH_s\cap UH_s^c $ occurs, there cannot exists a circuit in $A_{s,\f{C_1s}{2}}$  around $B_s$, that is $LH_s\cap UH_s^c \subset \cE^c$. This observation together with (\ref{eq:event_E_2}) and translation invariance yields
%the definition of $\cH_s(\al,\be)$
%
\begin{eqnarray}
\f{c_0}{4} > P\left( LH_s\cap UH_s^c \right)&\geq& P\left(LH_s\right)-P\left(UH_s\right)\nn\\
&=&P\left(\cH_{C_1s}(0,s)\right)-P\left(\cH_{C_1s}\left(s,\f{C_1s}{2}\right)\right),
\lab{eq:LH_UH_1}
\end{eqnarray}
%
%Observe that when $LH_s\cap UH_s^c $ occurs, there cannot exists a circuit in $A_{s,\f{C_1s}{2}}$  around $B_s$, that is $LH_s\cap UH_s^c \subset \cE^c$. That implies
%
%\begin{equation}
% P\left(LH_s\cap UH_s^c \right)\leq P\left(\cE^c\right).  
%\lab{eq:LH_UH_3}
%\end{equation}
%
%From (\ref{eq:event_E_2}) and (\ref{eq:LH_UH_1}) we obtain
%
%\begin{equation} 
%P\left(\cH_{C_1s}(0,s)\right)-P\left(\cH_{C_1s}(s,C_1s/2)\right)<c_{0}/4.
%\lab{eq:}
%\end{equation}
%
which contradicts (\ref{eq:prob_chi_with_constant}) and hence the assumption that $\al_t<s$ for all $t \in \big[4s,C_1s\big]$.  So there must exists some  $ s' \in [4s,C_1s]$, such that $\al_{s'}\geq s$. This proves the first assertion in Lemma~\ref{lem:condition_on_constants}.

%The second part of Lemma~\ref{lem:condition_on_constants} asserts the existence of $s_2$ with the desired properties. 
From the first part of this lemma there exists $s' \in [4s,C_1s]$ such that 
\begin{equation}
\al_{s'}\geq s\geq  \f{s'}{C_1}.
\lab{eq:s1_to_s1'}
\end{equation}
%
%where $C_1$ is as defined in (\ref{eq:condition_on_constants}). 
We shall prove the second part as well by contradiction. Suppose if possible $\al_t> 2\al_{2t/3}$ for all $t\geq s'$. By iterating this inequality we obtain
\begin{equation}
\al_{\left( \f{3}{2} \right)^k s'}>  2\al_{\left(\f{3}{2}\right)^{k-1} s'}> 2^k \al_{s'}\ge 2^k \f{s'}{C_1},
\lab{eq:inequality_on_alpha_k}
\end{equation}
 for all $k\geq 1$, where the last inequality follows from~(\ref{eq:s1_to_s1'}). Since $\al_s\leq \f{s}{4}$ for all $s\geq 1$, we have for all $ k\geq 1$ that
\begin{equation}
\al_{\left(\f{3}{2}\right)^k s'}\leq \f{1}{4}\left(\f{3}{2}\right)^k s'.
\lab{eq:property_of_alpha_k}
\end{equation}
The inequalities (\ref{eq:inequality_on_alpha_k}) and (\ref{eq:property_of_alpha_k}) implies that for all $k\geq 1$,
\begin{equation}
\f{1}{4}\left(\f{3}{2}\right)^k > \f{2^k}{C_1},
\lab{eq:lowerbound_C1}
\end{equation}
which contradicts the fact that $C_1< \infty$. Hence the statement that $\al_t> 2\al_{2t/3}$ for all $t\geq s'$ is not true. In particular  $\al_{\left(\f{3}{2}\right)^k s'}>  2\al_{\left(\f{3}{2}\right)^{k-1} s'}$ is not true for all $k \geq 1$.

Let $k^*:=\min\big\{k \in \bN: \al_{\left(\f{3}{2}\right)^k s'}\leq  2\al_{\left(\f{3}{2}\right)^{k-1} s'}\big\}$. By definition of $k^*$ 
\begin{equation}
\al_{\left(\f{3}{2}\right)^{k^{\star}} s'}\leq 2\al_{\left(\f{3}{2}\right)^{{k^{\star}}-1} s'}.
\end{equation}
Again by definition of $k^*$ and the argument leading to (\ref{eq:lowerbound_C1}) we have $\left(\f{4}{3}\right)^{k^{\star}-1}<\f{C_1}{4}$, which implies that $k^{\star}\leq \lfloor \log_{\f{4}{3}}{\f{C_1}{4}} \rfloor +1$. 

Let $s_2:=\left(\f{3}{2}\right)^{k^{\star}}s'$. Observe that $s_2\geq s'$ and $s_2 \leq \left(\f{3}{2}\right)^{  \lfloor \log_{\f{4}{3}}{\f{C_1}{4}} \rfloor +1} s'$.
%\leq (C_1/3)^{\log_{\f{4}{3}}{\f{3}{2}}} s_1' $. 
Let $C_1' = \left(\f{3}{2}\right)^{  \lfloor \log_{\f{4}{3}}{\f{C_1}{4}} \rfloor +1}$. Thus we have found $s_2 \in [s', C_1' s']$ such that $\al_{s_2}\leq 2\al_{\left(\f{2}{3}\right)s_2}$.

This proves second part of Lemma~\ref{lem:condition_on_constants} and completes the proof of  Proposition \ref{prop:sequence_of_good_scales}.\qed
\begin{figure}[h!]
  \centering
\hspace{-0.4in}
\includegraphics[width=0.85\linewidth]{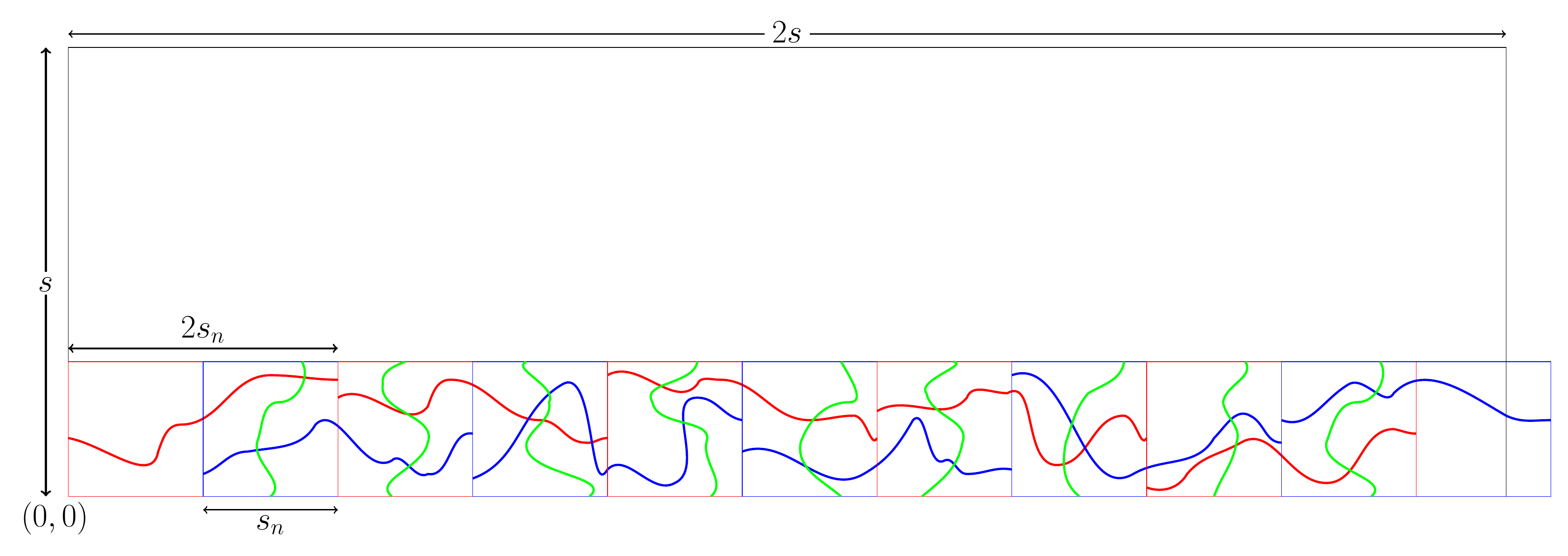}
  \captionof{figure}{The event $F_s$, occurrence of which gives a left-right crossing in $[0,2s] \times [0,s]$.}
  \label{fig:overlapping_paths_2}
\end{figure}

{\bf Proof of Proposition \ref{prop:general_s_2}. } Let $\{s_n\}_{n\geq 1}$ be the sequence of scales as in Proposition \ref{prop:sequence_of_good_scales}. For any $s\geq 1$ let $k = k(s)$ be such that $s_k \leq s<s_{k+1}$. Since $s_{k+1} \leq C s_k$ with $C$ as in Proposition~\ref{prop:sequence_of_good_scales} we have
\[ [0,2s]\times[0,s_k]\subset \bigcup\limits_{i=0}^{n_1}\left((is_k,0)+[0,2s_k]\times[0,s_k]\right), \]
where $n_1=\lfloor2s/s_k\rfloor+1\leq 2 C +1$.

Let $F_s$ be the event (see Figure ~\ref{fig:overlapping_paths_2}) that there is a left to right crossing in each of the rectangles $(is_k,0)+[0,2s_k]\times[0,s_k]$ for $i=0,1,\cdots,n_1$ and top to down crossing in each of the squares $(is_k,0)+[0,s_k]\times[0,s_k]$ for $i=1,2,\cdots,n_1$. Clearly $F_s \subset LR_s(2)$ and by Proposition \ref{prop:sequence_of_good_scales} we have $C_{s_k}(2) \geq  P(\cA_{s_k}) \geq c_2$. It follows by the FKG inequality that $C_s(2)\geq c_2^{n_1+1}c_0^{n_1} > 0.$ This proves Proposition~\ref{prop:general_s_2}.  \qed

\subsection{ Proof of Theorem~\ref{thm:RSW_nonpercolation}~(\ref{rsw_lim}) for eRCM } As in the proof of the first part it suffices to show the result for $\rho = 2$. We first complete the proof using the following lemma which will be proved subsequently using techniques similar to that used to prove Lemma~\ref{lem:condition_on_constants}.
\begin{lem}
Suppose that $\lim\limits_{s\rar \infty}C_s(1)=1$. Then for any fixed $\ep>0$ there exists $\eta \in (0,\f{1}{4})$ such that for all $s$ sufficiently large we have
\begin{equation}
P\left[\mbox{there exists a circuit around } B_{\eta s}\mbox{ in the annulus } A_{\eta s, \f{s}{4}}\right]>1-\ep.
\lab{eq:circuit_eta}
\end{equation}
\lab{lem:crossing_eta_annulus}
\end{lem}
\begin{figure}[h!]
  \centering
\hspace{-0.4in}
\includegraphics[width=0.58\linewidth]{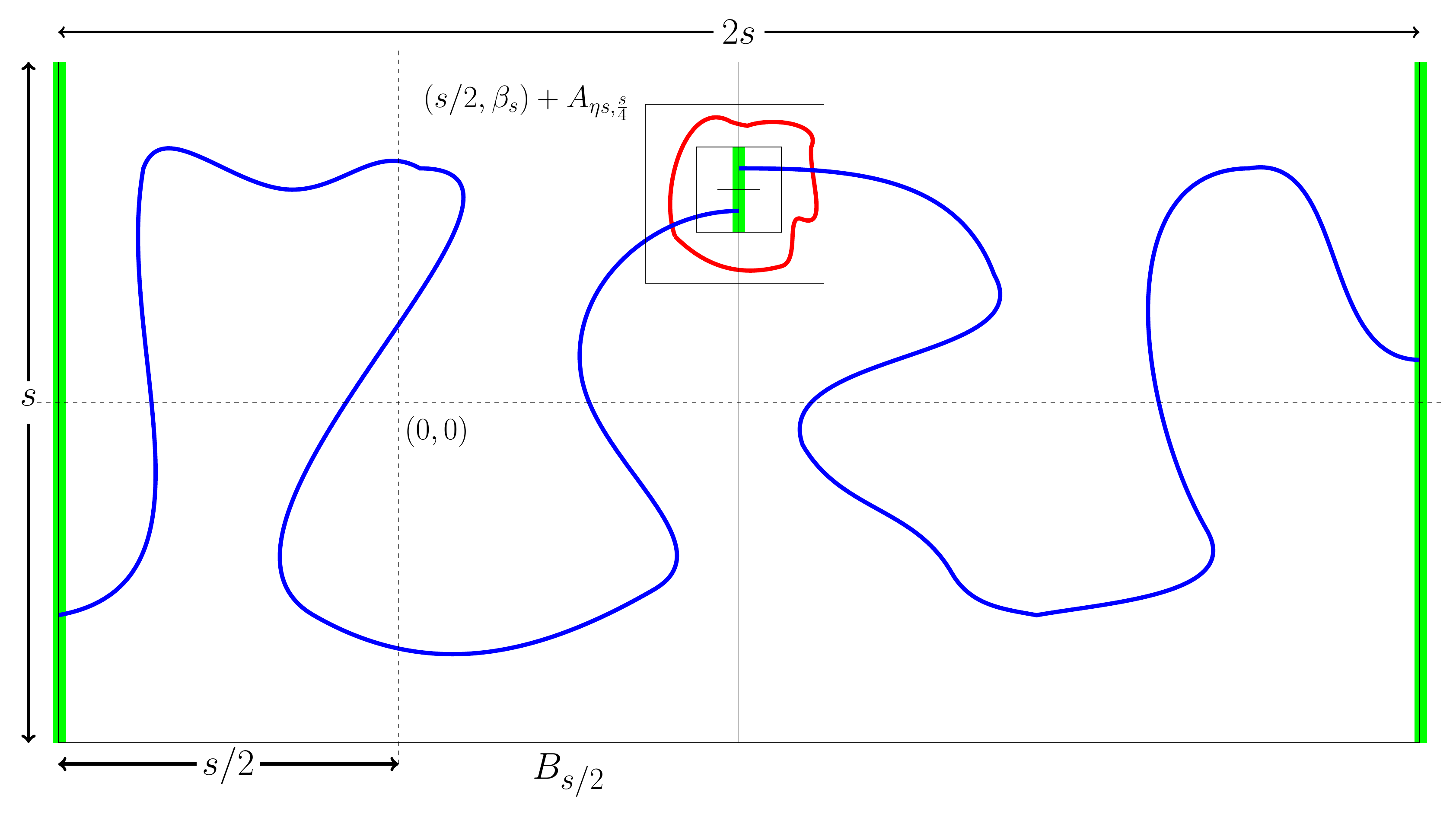}
  \captionof{figure}{A realisation of the event $\cT_s$}
  \label{fig:event_Ts}
\end{figure}
Let $\eta>0$ be as in Lemma~\ref{lem:crossing_eta_annulus}. Divide the side $\{\f{s}{2}\}\times [-\f{s}{2},\f{s}{2}]$ into intervals labeled $J_i(\eta, s)$, $i=1, 2, \ldots , \kappa$, $\kappa= \lfloor \f{1}{2\eta}\rfloor + 1$, of length $2\eta s$ (except for one interval that is of length at most $2 \eta s$). For $i=1, 2, \ldots , \kappa,$ let $\cH(J_i(\eta , s))$ be the event that there exists a path in the box $ B_{\f{s}{2}}$ from left to $ J_i(\eta, s)$. Clearly  $LR_s(1)= \bigcup_{i=1}^{\kappa} \cH(J_i(\eta, s)) $.  Using square root trick there exists $\be_{s}\in \left[-\f{s}{2},\f{s}{2}\right]$ satisfying $-\left(\f{1}{2}-\eta\right)s\leq \be_{s}\leq \left(\f{1}{2}-\eta\right)s$ such that,
\begin{eqnarray}
P\left( \cH_s(\be_{s}-\eta s,\be_{s}+\eta s )\right) = \max_{i} P\left(\cH(J_i(\eta , s))\right) & \geq & 1 - \left(1- P\left(\bigcup_{i=1}^{\kappa} \cH(J_i(\eta, s))\right) \right)^{\f{1}{\kappa}} \nn\\
& = & 1 - \left(1-C_s(1) \right)^{\f{1}{\kappa}}.
\lab{eq:crossing_to_interval_at_beta_s}
\end{eqnarray}
Let $R_s(\eta)$, $\cA'_{\eta s, \f{s}{4}}$ be the events that there exists a path from $\{\f{s}{2}\}\times [\beta_s-\eta s,\beta_s+\eta s]$ to right  in $(s,0)+B_{\f{s}{2}}$ and there exists a circuit in $(\f{s}{2}, \beta_s)+A_{\eta s, \f{s}{4}}$ respectively (see Figure~\ref{fig:event_Ts}). By translation and rotation invariance $P\left( R_s(\eta)\right)=P\left( \cH_s(\beta_{s}-\eta s,\beta_{s}+\eta s )\right)$.
For any $\ep > 0$, by Lemma~\ref{lem:crossing_eta_annulus} we have $P(\cA'_{\eta s, \f{s}{4}}) > 1 - \ep$ for all $s$ sufficiently large. Let $\cT_s:= \cH_s(\beta_{s}-\eta s,\beta_{s}+\eta s)\cap R_s(\eta)\cap \cA'_{\eta s, s/4}$. By the FKG inequality and (\ref{eq:crossing_to_interval_at_beta_s}) we obtain
\begin{eqnarray*} 
C_s(2) & \geq & P(\cT_s) = P\left(\cH_s(\beta_{s}-\eta s,\beta_{s}+\eta s)\cap R_s(\eta)\cap \cA'_{\eta s, s/4}\right) \\ 
 & \geq &  \left[1-(1-C_s(1))^{\f{1}{\kappa}} \right]^2 (1-\ep) \to (1-\ep),
\end{eqnarray*}
as $s \to \infty$. The result now follows since $\ep > 0$ is arbitrary. \qed

\subsection{Proof of Lemma~\ref{lem:crossing_eta_annulus}}

 Fix $\ep, \tau \in (0, 1)$. Since $\lim\limits_{s \to \infty}C_s(1) = 1$, there exists a $s_0 > 0$ such that $\inf\limits_{s \geq s_0} C_s(1) > 0$ and hence by Theorem~\ref{thm:RSW_nonpercolation}~(\ref{rsw_inf}) $\inf\limits_{s \geq s_0} C_s(4) > 0$. By the FKG inequality, 
\begin{equation}
c := \inf\limits_{s \geq s_0} P\left(\cA_ s\right) \geq \inf\limits_{s \geq s_0} \left(C_ s(4)\right)^4 > 0.
\lab{eq:prob_scaled_annulus_eta}
\end{equation}
Choose $\eta \in (0, \f{1}{4})$ satisfying 
\begin{equation}
\left(1-\f{c}{2}\right)^{\lfloor\log_5\f{1}{8\eta}\rfloor} < \ep.
\lab{eq:condition_on_constants_eta}
\end{equation}
From (\ref{eq:prob_scaled_annulus_eta}) there is a $s_1 > s_0$ such that for all $s\geq s_1$ and $i=1,2,\cdots, \lfloor\log_5 \f{1}{8\eta}\rfloor$ we have
\begin{equation}
 P\left(\cA_{\f{5^i}{2}\eta s}\right) \geq  c.
\lab{eq:prob_scaled_annulus_i}
\end{equation}
Using Corollary~\ref{cor:edge_annulus} choose $s_2 > s_1$ such that for all $s\geq s_2$ we have
\begin{equation}
 P\left(L_{\f{5^i}{4}\eta s}(\tau) \right)\leq \f{c}{2}.
\lab{eq:measurable_L_eta}
\end{equation}
Combining (\ref{eq:prob_scaled_annulus_i}) and (\ref{eq:measurable_L_eta}) we can write for  $i=1,2,\cdots, \lfloor\log_5 \f{1}{8\eta}\rfloor$,
\begin{equation}
 P\left(\cA_{ \f{5^i}{2}\eta s}\cap L^c_{\f{5^i}{4}\eta s}(\tau)\right) \geq \f{c}{2}.
\lab{eq:combined_A_L_eta}
\end{equation}
Let $\cE_s$ be the event that there exists a circuit around $B_{\eta s}$ in the annulus $A_{\eta s, \f{s}{4}}$. Observe that if $\cA_{ \f{5^i}{2}\eta s}\cap L^c_{\f{5^i}{4}\eta s}(\tau)$ occurs for some $ i=1,2,\cdots, \lfloor\log_5 \f{1}{8\eta}\rfloor$,  then $\cE_s$ will also occur. By the same argument as in Lemma~\ref{lem:condition_on_constants} appearing below (\ref{eq:event_E_0}), the events  $\cA_{ \f{5^i}{2}\eta s}\cap L^c_{\f{5^i}{4}\eta s}(\tau)$, $ i=1,2,\cdots, \lfloor\log_5 \f{1}{8\eta}\rfloor$ are independent. Using the above two observations along with (\ref{eq:combined_A_L_eta}) and (\ref{eq:condition_on_constants_eta}) yields
\begin{eqnarray*}
P(\cE_s^c) & \leq & 
P\left(\bigcap_{i=1}^{\lfloor\log_5\f{1}{8\eta}\rfloor}  \left(\cA_{ \f{5^i}{2}\eta s}\cap L^c_{\f{5^i}{4}\eta s}(\tau)\right)^c\right) \\
& = & \prod\limits_{i=1}^{\lfloor\log_5 \f{1}{8\eta }\rfloor} P\left( \left(\cA_{ \f{5^i}{2}\eta s}\cap L^c_{\f{5^i}{4}\eta s}(\tau)\right)^{c}\right)< \left(1-\f{c}{2}\right)^{\lfloor\log_5\f{1}{8\eta}\rfloor} <\ep.
\hfill \qed
\end{eqnarray*}
\subsection{Proof of Theorem~\ref{thm:RSW_nonpercolation}~(\ref{non_perc})}
The proof follows by a renormalization argument that uses the RSW Lemma (Theorem~\ref{thm:RSW_nonpercolation}~(\ref{rsw_lim})) and Proposition~\ref{prop:largest_edge_length_in_box} on the length of the longest edge in the graph $G_{\lam}^e$. Suppose $\th^e(\lam)>0$, that is, the graph $G_\lam^e$ percolates. For $u < n$ and fixed $\tau \in \left(\f{2}{c-2}, 1\right)$  define the events 
\begin{equation}
E(u, n) : =\Big\{B_{u}\mathrel{\mathop{\llra}^{(2n,0)+B_{4n}}} (4n,0)+B_{u}\Big\}\mbox{ and } \tilde{E}(u, n) : = E(u, n) \cap \cA_u\cap \tilde{\cA}_u\cap \{M_{4n}\leq  n^\tau\},
\lab{eq:edge_event}
\end{equation}
where $\cA_u$ $(\tilde{\cA}_u \mbox{ resp.})$ be the event that there exists a circuit around $B_u$ $((4n,0)+B_u \mbox{ resp.})$ in the annulus $A_{u,2u}$ $((4n,0)+A_{u,2u} \mbox{ resp.})$, $M_{4n}:=$ length of the longest edge intersecting the box $(2n,0)+B_{4n}$. We complete the proof using the following proposition, the proof of which shall be provided subsequently.
\begin{prop}
For the random graph $G_\lam^e$ with the connection function $g$ satisfying $g(r)=O(r^{-c})$ as $r \to \infty$ for $c>4$, if $\th^e(\lam)>0$ then there exists a sequence $\{u_n\}_{n\geq 1}$ satisfying $u_n\rar \infty$  as $n \rar \infty$ and $u_n < \f{n}{2}$ such that
 \begin{equation}
\lim_{n\rar \infty}P\left( E(u_{2n}, n) \right)=1.
\lab{eq:box_to_shifted_box}
\end{equation}
\lab{prop:box_to_shifted_box}
\end{prop}
Let $u_n$ be as in Proposition~\ref{prop:box_to_shifted_box}. We define a coupled nearest neighbor bond percolation model on $4n\bZ^2$. The edge $((0,0),(4n,0))$ is said to be open if $\tilde{E}(u_{2n}, n)$ occurs in $G_\lam^e$. An edge that is not open is designated closed. For any two nearest neighbors $z_1, z_2\in \bZ^2$ we can define an open edge between $4nz_1,4nz_2$ in an analogous manner. Otherwise the edge is said to be closed. Denote by $\tilde{G}_n$ the graph on the vertex set $4n\bZ^2$  formed by the open edges. By translation and rotational invariance, any edge is open has probability $P\left(\tilde{E}(u_{2n}, n)\right)$. The edge $((0,0),(4n,0))$ being open does not depend on the configuration of points of $\cP_\lam$ outside $(2n,0) + B_{5n}$. Thus the status of the edge $((0,0),(4n,0))$ can influence that of at most forty neighboring edges. By Theorem $0.0$ by Liggett et.al \citep{Liggett1997} for finitely dependent nearest neighbor bond percolation model on $4n\bZ^2$, $n \in \bN$, there exists a constant $q_0\in(0,1)$ such that the random graph $\tilde{G}_n$ percolates whenever, 
\begin{equation}
P\left(\tilde{E}(u_{2n}, n)\right)>q_0.
\lab{eq:P_tildeE}
\end{equation} 
By the FKG inequality $P\left( \cA_s \right) \geq C_s(4)^4$. Hence by translational invariance, Theorem~\ref{thm:RSW_nonpercolation}\,$ (\ref{rsw_lim})$, Proposition~\ref{prop:box_to_shifted_box} and  Proposition~\ref{prop:largest_edge_length_in_box} we have
%
%\begin{equation}
%\min\left\{ P\left(E(u_{2n}, n)\cap \left\{ M_{4n} \leq n^\tau \right\} \right), P\left(\cA_{u_{2n}}\cap \left\{ M_{4n} \leq n^\tau \right\} \right), P\left( \tilde{\cA}_{u_{2n}}\cap \left\{ M_{4n} \leq n^\tau \right\}\right) \right\} \to 1,
%\lab{eq:E_un_script_A_script_A_tilde}
%\end{equation}
%
%as $n \to \infty$. Since $E(u_{2n}, n),  \cA_{u_{2n}}, \tilde{\cA}_{u_{2n}}$  are increasing events, by the FKG inequality 
%
\begin{equation}
P\left(\tilde{E}(u_{2n}, n)\right) \to 1, \qquad \mbox{as } n \to \infty.
\lab{eq:function_geq_q0}
\end{equation}
%
%\begin{eqnarray}
%P\left(\tilde{E}(u_{2n}, n)\right) & = & P\left(E(u_{2n}, n)\cap \cA_{u_{2n}} \cap \tilde{\cA}_{u_{2n}} \cap \left\{ M_{4n} \leq n^\tau \right\}\right) \nn \\
%& \geq &  P\left(E(u_{2n}, n)\cap \left\{ M_{4n} \leq n^\tau \right\} \right) P\left(\cA_{u_{2n}} \cap \left\{ M_{4n} \leq n^\tau \right\} \right)P\left( \tilde{\cA}_{u_{2n}}\cap \left\{ M_{4n} \leq n^\tau \right\}\right) \nn \\
% & \to & 1, \qquad \mbox{as } n \to \infty.
%\lab{eq:function_geq_q0}
%\end{eqnarray}
%
For $n\in \bN$ and $\lam > 0$ define $f_{n}(\lam) := P\left(\tilde{E}(u_{2n}, n)\right)$. Let $X_h$ be a Poisson random variable with mean $100n^2h$. Then a simple coupling argument shows that $\big|f_n(\lam+h)-f_n(\lam)\big|\leq P\left(X_h \geq 1\right) \to 0$ as $h \to 0$. So $f_n$ is continuous. 

Let $\Lam:=\{\lam: \th^e(\lam)>0\}$ be the set of parameters $\lam > 0$ for which $G^e_\lam$ percolates. Since $G^e_\lam$ percolates if $\tilde{G}_n$ does we have from (\ref{eq:P_tildeE}) that $\bigcup_n f_n^{-1}(q_0,1] \subset \Lam$. On the other hand if $\lam\in \Lam$, then by (\ref{eq:function_geq_q0}) there exists an $n_0\in \bN$ such that $\lam \in f^{-1}_{n_0}(q_0,1]$ and hence $\Lam \subset \bigcup_n f_n^{-1}(q_0,1]$. It follows that $\Lam=\bigcup_n f_n^{-1}(q_0,1]$. Since the functions $f_n$ are continuous, $\Lam$ is an open set. This completes the proof of Theorem~\ref{thm:RSW_nonpercolation}~(\ref{non_perc}). It remains to prove the Proposition~\ref{prop:box_to_shifted_box}. \qed

{\bf Proof of Proposition~\ref{prop:box_to_shifted_box}. } Since $\th^e(\lam)>0$, there exists almost surely an infinite connected component in $G_\lam^e$. Hence for any sequence $\{u_n\}_{n\geq1}$ satisfying $u_n\rar \infty$  as $n\rar \infty$ such that
\begin{equation}
P\left(\cC\;  \mbox{intersects} \; B_{u_n}\right)\rar 1 \mbox{ as } n\rar \infty,
\lab{eq:inf_comp_intersects}
\end{equation}
where $\cC$ is an infinite connected component in $G_\lam^e$. Fix one such sequence for which $u_n<\f{n}{2}$. An immediate consequence of (\ref{eq:inf_comp_intersects}) is that
\begin{equation}
P\left( B_{u_n} \llra \partial B_n\right)\rar 1 \mbox{ as } n\rar \infty.
\lab{eq:crossing_the_boundary}
\end{equation}
%
%Recall that $M_n$ is the length of the longest edge in $G_\lam$ \red{that intersects} $B_n$. By Proposition~\ref{lem:largest_edge_length_in_box} and the FKG inequality we have
%
%\begin{align}
%P\left(B_{u_n}\mathrel{\mathop{\llra}^{B_n}} \partial B_n \right)\to 1.
%\lab{eq:null_prob}
%\end{align}
%
Let $v_n, w_n$ be as in Lemma~\ref{lem:hit_the_interval} below. Define the events
\begin{equation}
H_n:=B_{u_{2n}}\mathrel{\mathop{\llra}^{B_{2n}}}\{2n\}\times \big[v_{2n},  w_{2n}\big],\,\,
H'_n:=(4n, 0)+B_{u_{2n}}\mathrel{\mathop{\llra}^{(4n,0)+B_{2n}}} \{2n\} \times\big[v_{2n}, w_{2n}\big] .\nn
\end{equation}
For $n \in \bN$ let $\hat{\cA}_n$ be the event that there exists a circuit around $\left(2n,\f{1}{2}(v_{2n}+w_{2n})\right)+B_{\f{1}{2}(w_{2n}-v_{2n})}$ within the annulus $\left(2n,\f{1}{2}(v_{2n}+w_{2n})\right) + A_{\f{1}{2}(w_{2n}-v_{2n}), (w_{2n}-v_{2n})}$. We denote the annulus $\left(2n,\f{1}{2}(v_{2n}+w_{2n})\right) + A_{\f{1}{2}(w_{2n}-v_{2n}), (w_{2n}-v_{2n})}$ by $A(n, v_{2n}, w_{2n})$. By the Lemma~\ref{lem:hit_the_interval} as $n \rar \infty$
\begin{equation}
P\left(H_n \right)=P\left(H'_n \right)\rar 1 \qquad \mbox{and} \qquad A(n, v_{2n}, w_{2n})  \subset (2n,0)+B_{2n}.
\lab{eq:two_paths_1}
\end{equation}
By translation invariance, the FKG inequality and Theorem~\ref{thm:RSW_nonpercolation}~(\ref{rsw_lim}) we have
\begin{equation}
\lim_{n \rar \infty}P\left(\hat{\cA}_n\right)=1.
\lab{eq:circuit_1}
\end{equation}
\begin{figure}
  \centering
\includegraphics[width=0.6\linewidth]{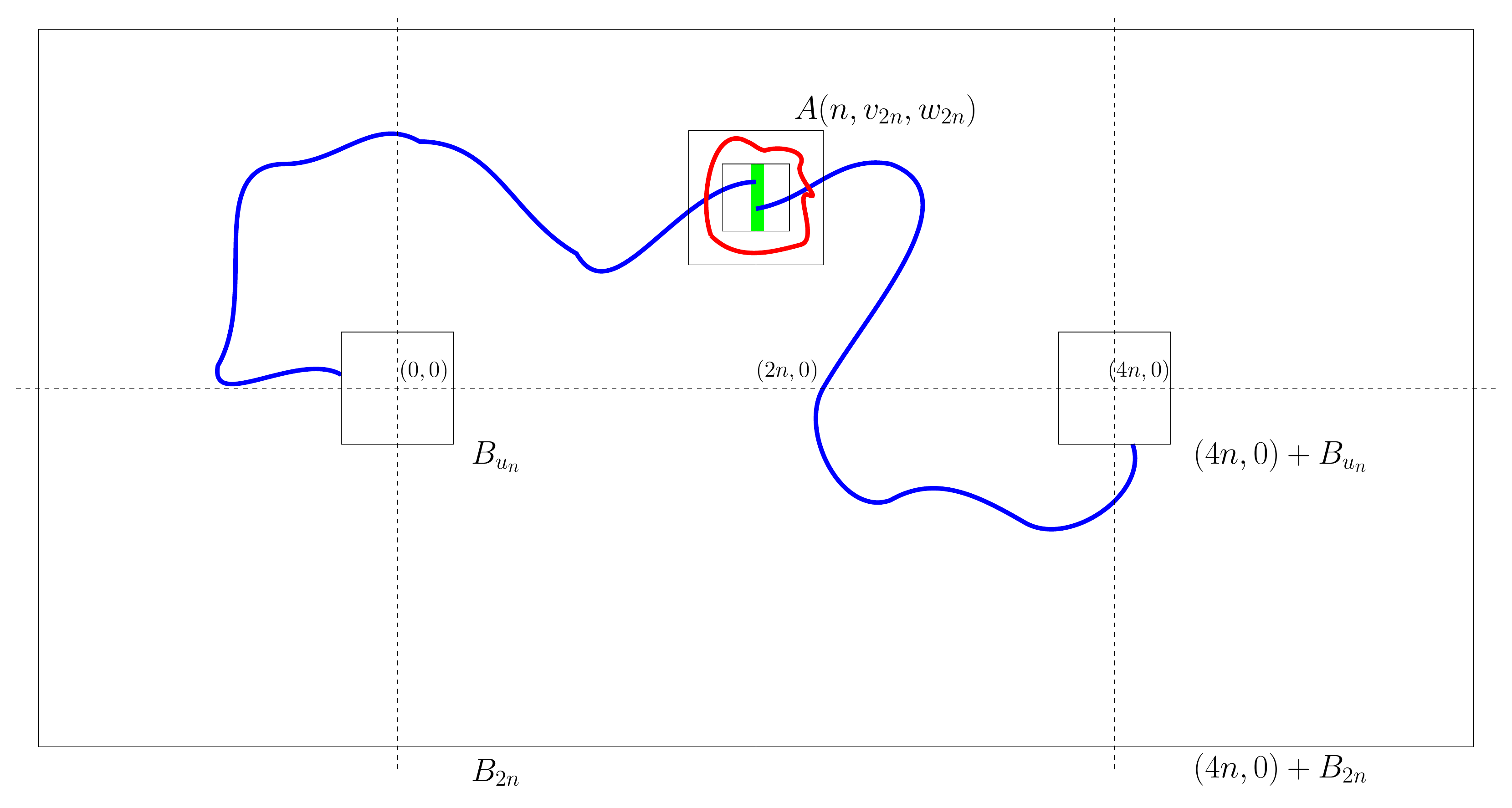}
  \captionof{figure}{A realisation of the event $H_n \cap H'_n \cap \hat{\cA}_n$}
  \label{fig:event_En}
\end{figure}
Observe that the paths that enable the events $H_n, H'_n$ must intersect the circuit in $A(n, v_{2n}, w_{2n})$ (see Figure~\ref{fig:event_En}) and hence
\begin{equation}
H_n \cap H'_n \cap \hat{\cA}_n \subset \big\{B_{u_{2n}}\mathrel{\mathop{\llra}^{(2n,0)+B_{4n}}} (4n,0)+B_{u_{2n}}\big\}.
\lab{eq:circuit_path}
\end{equation}
Proposition~\ref{prop:box_to_shifted_box} now follows from (\ref{eq:two_paths_1})--(\ref{eq:circuit_path}). \qed
\begin{lem}
 Consider the random graph $G_{\lam}^e$ with the connection function $g$. If $\th^e(\lam)>0$ the for any $k\in \bN$ there exists sequence $\{v_{n}\}_{n\geq 1}, \{w_{n}\}_{n\geq 1}$  satisfying $\left[ v_{2n} - \f{1}{2}(w_{2n}-v_{2n}),  w_{2n} +\f{1}{2}(w_{2n}-v_{2n}) \right] \subset [- 2n, 2n]$ such that the following holds as $n\rar \infty$
\begin{equation}
P\left(B_{u_{2n}}\mathrel{\mathop{\llra}^{B_{2n}}}\{2n\}\times \big[v_{2n}, w_{2n}\big] \right) \rar 1.
\lab{eq:half}
\end{equation}
\lab{lem:hit_the_interval}
\end{lem}
{\bf Proof of Lemma~\ref{lem:hit_the_interval}.} Consider the square $B_{2n}$. By rotational invariance the probability of having a path from $B_{u_{2n}}$ to any of the eight half intervals on the sides of the square $B_{2n}$ are same. In other words
\begin{equation}
P\left(B_{u_{2n}}\mathrel{\mathop{\llra}^{B_{2n}}} A_i\right)=P\left(B_{u_{2n}}\mathrel{\mathop{\llra}^{B_{2n}}} A_j\right)
\lab{eq:equal_hands}
\end{equation}
for $i,j\in \{1,2,\cdots, 8\}$, where $A_1=\{2n\}\times [0,2n],A_2=[0,2n]\times\{2n\},\ldots ,A_8=\{2n\}\times [-2n,0]$ are the half intervals on the sides of $B_{2n}$. Applying square root trick and using (\ref{eq:crossing_the_boundary}) we have
\begin{equation}
P\left(B_{u_{2n}}\mathrel{\mathop{\llra}^{B_{2n}}} A_i \right)\geq 1-\left(1-P(B_{u_{2n}}  \llra \partial B_{2n} )\right)^{\f{1}{8}} \to 1, 
\lab{eq:lb_equal_hand_prob}
\end{equation}
as $n \to \infty$, for all $i=1, \ldots ,8$.
For any $n\in \bN$ and  $\th\in[0,2n]$
\begin{equation}
\big\{B_{u_{2n}}\mathrel{\mathop{\llra}^{B_{2n}}}\{2n\}\times[0,\th] \big\} \cup \big\{B_{u_{2n}}\mathrel{\mathop{\llra}^{B_{2n}}}\{2n\}\times[\th,2n] \big\}=\left\{B_{u_{2n}}\mathrel{\mathop{\llra}^{B_{2n}}} A_1\right\}.\nn
\end{equation} 
For $\th\in [0,2n]$ define the functions
\begin{equation}
L_{2,n}(\th):=P\left(\left\{B_{u_{2n}}\mathrel{\mathop{\llra}^{B_{2n}}}\{2n\}\times[0,\th] \right\} \right),\, U_{2,n}(\th):=P\left(\left\{B_{u_{2n}}\mathrel{\mathop{\llra}^{B_{2n}}}\{2n\}\times[\th,2n] \right\} \right).\nn
\end{equation}
%
%From (\ref{eq:lb_equal_hand_prob}) and another application of square root trick we obtain
%
%\begin{equation}
%\max\big\{L_{2,n}(\th), U_{2,n}(\th)\big\} \geq     1-\left(1-P\left(B_{u_{2n}}\mathrel{\mathop{\llra}^{B_{2n}}} A_1 \right)\right)^{\f{1}{2}} \to 1.
%\lab{eq:max_Un_Ln_limit}
%\end{equation}
%
%From (\ref{eq:lb_equal_hand_prob}) and (\ref{eq:max_Un_Ln}) we obtain, 
%
%\begin{equation}
%\lim_{n\rar \infty}\max\big\{L_{2,n}(\th), U_{2,n}(\th)\big\} =1.
%\lab{eq:max_Un_Ln_limit}
%\end{equation}
%
Observe that $U_{2,n}(0)>L_{2,n}(0)=0,\,L_{2,n}(2n)>U_{2, n}(2n)=0,\, L_{2, n}$ is non-decreasing and $U_{2,n}$ is non-increasing and $L_{2, n},\,U_{2, n}$ are continuous.
By the properties of $L_{2, n},\,U_{2, n}$  there exists $\al_{2n} \in (0,2n)$ such that  $U_{2,n}(\th)\geq L_{2, n}(\th)$ for $\th<\al_{2n},\,L_{2,n}(\th)\geq U_{2,n}(\th)$ for $\th>\al_{2n}$. Further from (\ref{eq:lb_equal_hand_prob}) and another application of square root trick we obtain 
\begin{equation}
U_{2,n}(\al_{2n})=L_{2,n}(\al_{2n}) \to 1.
\lab{eq:Max_U_L}
\end{equation}
%
%It follows from (\ref{eq:Max_U_L}) and (\ref{eq:max_Un_Ln_limit}) that
%
%\begin{equation}
%\lim_{n\rar\infty}U_{2, n}(\al_{2n}) =  \lim_{n\rar\infty}L_{2, n}(\al_{2n})=1.
%\lab{eq:limit_Un}
%\end{equation}
%
%It follows from \ref{eq:Max_U_L} that 
%
%\begin{equation}
%P\left(B_{u_{2n}}\mathrel{\mathop{\llra}^{B_{2n}}}\{2n\}\times [\al_{2n},2n]\right) = P\left(B_{u_{2n}}\mathrel{\mathop{\llra}^{B_{2n}}}\{2n\}\times [0, \al_{2n}]\right) \rar 1,
%\lab{eq:upper_part}
%\end{equation}
%
Set $\be_{2n}= 2n-\al_{2n} > 0$ and let $\gamma_{2n}:= \f{\al_{2n} \wedge \be_{2n}}{16}$. Let $\{2n\}\times [0, \al_{2n}] = \bigcup_{j=1}^k I_j^{(n)}$ where $I_j^{(n)}= \{2n\}\times [(j-1)\gamma_{2n}, j\gamma_{2n}]$, for $j=1,2, \ldots ,(k-1)$ and $I_k^{(n)}= \{2n\}\times [(k-1)\gamma_{2n}, k \gamma_{2n}]$, $k=\lfloor \f{\al_{2n}}{\gamma_{2n}}\rfloor +1 $ (\ref{eq:Max_U_L}) together with an application of square root trick yields
\begin{equation}
\max_{j\in [k]} \Big\{P\left(B_{u_{2n}}\mathrel{\mathop{\llra}^{B_{2n}}}I_j^{(n)} \right) \Big\}\to 1.
\lab{eq:two_quarter_interval}
\end{equation}
%
%%%%%%%%%%%% removed
\remove{
\begin{equation}
\max \Big\{P\left(B_{u_{2n}}\mathrel{\mathop{\llra}^{B_{2n}}}\{2n\}\times \big[\al_{2n}, \f{3}{4}\al_{2n}\big] \right), P\left(B_{u_{2n}}\mathrel{\mathop{\llra}^{B_{2n}}}\{2n\}\times \big[ \f{3}{4}\al_{2n}, 2n\big]\right) \Big\}\to 1.
\lab{eq:two_quarter_interval}
\end{equation}
}
%%%%%%%%%%%% removed
%
as $n\to \infty$. Let $t(n) = \argmax_{j\in [k]} P\left(B_{u_{2n}}\stackrel{B_{2n}}{\llra} I_j^{(n)} \right)$. (\ref{eq:half}) follows by taking $[v_{2n},w_{2n}] = I_{t(n)}^{(n)}.$ \qed

\subsection{Proof of Theorem~\ref{thm:PhaseTransition}~(\ref{np_iercm})}

The proof of non-trivial phase transition for ieRCM is similar to that of eRCM. The analysis is bit more involved due to the presence of weights at the vertices. From Theorem 3.2 $(a2)$ in \cite{Deprez2018} with $d=2$ we have that a phase transition occurs in $H_\lam$ provided $\al > 2, \al \be > 4$. However, this condition is required only to show that $\tilde{\lam}_c > 0$ while $\tilde{\lam}_c < \infty$ holds over the entire parameter space. Clearly, $H^e_\lam$ percolates  if $H_\lam$ does and hence $\tilde{\lam}_c^e < \infty$. 

We now show $\tilde{\lam}^e_c> 0$ using the same technique as in the proof of Theorem~\ref{thm:PhaseTransition}~(\ref{np_ercm}). With the same notations as in that proof, we start with the inequality analogous to (\ref{eq:theta_tilde_2}). 
\begin{eqnarray}
\tilde{\th}^e(\lam) 
%
%&\leq &  P^o(\mbox{ there is a self-avoiding  path on } n \mbox{ vertices in } H_{\lam}^e)\nn\\
%
%& \leq &  E^o \left[\sum\limits_{\textbf{x} \in \cP_{\lam,\neq}^n} 1\{x_0 \to x_1 \to x_2 \to \cdots \to x_n \mbox{ occurs}\} \right] \nn\\
%
& \leq &   \sum\limits_{k=0}^{n-1}  \sum\limits_{B \in \cB_k}  E^o \left[   \sum\limits_{\textbf{x} \in \cP_{\lam,\neq}^n} 1\{B(\textbf{x}) = B\}\right].
\lab{eq:thetaTilde_1}
\end{eqnarray}
%
%%%%%%%%%removed 
%
\remove{ For a block of the form $B_j := \{i: i_{j-1}\leq i \leq i_j\}$ of size larger than two, let 
\[ A_j = \{\textbf{y}=(y_1,  \ldots , y_n): \overline{y_{i_{j-1} + 2r - 2}y_{i_{j-1} + 2r - 1}} \mbox{ intersects } \overline{y_{i_{j-1} + 2r}y_{i_{j-1} + 2r + 1}} \mbox{ for all } r = 1, \ldots , (i_j - i_{j-1} - 1)/2 \}. \]
The intersection above is understood to be at an interior point of the line segments. For a block of size two we set $1_{A_j} \equiv 1$. }
%
%%%%%%%%%removed 
%
Let $\sigma(W):=\sigma\{W_x : x\in \cP_\lam\}$ be the sigma algebra generated by the weights. Conditioning on $\cP_{\lam}$, $\sigma(W)$ and then using (\ref{eq:connection_function}) along with the inequality $1 - e^{-x} \leq x \wedge 1$ for $x \geq 0$ we can bound the expectation inside the sum in (\ref{eq:thetaTilde_1}) as follows. 
\begin{eqnarray}
\lefteqn{E^o \left[   \sum\limits_{\textbf{x} \in \cP_{\lam,\neq}^n} 1\{B(\textbf{x}) = B\}\right]
 =    E^o \left[   \sum\limits_{\textbf{x} \in \cP_{\lam,\neq}^n} E^o \left[\prod\limits_{j=1}^{k+1} 1\{B_j(\textbf{x}) = B_j\}\bigg\vert \cP_{\lam}, \sigma(W)\right] \right] } \nn\\
 & = &   E^o \left[   \sum\limits_{\textbf{x} \in \cP_{\lam,\neq}^n} E^o \left[\prod\limits_{j=1}^{k+1} \prod\limits_{l\in B_j^e}g( x_{l-1}, x_l) 1_{A_j}(\textbf{x})\bigg\vert \cP_\lam\right] \right]  \nn\\
 & \leq &  E^o \left[   \sum\limits_{\textbf{x} \in \cP_{\lam,\neq}^n}   E^o\left[\prod\limits_{j=1}^{k+1} \prod\limits_{l\in B_j^e}  \left( \f{\eta  W_{x_{l-1}}W_{x_l}}{|x_l-x_{l-1}|^\al} \wedge 1\right) 1_{A_j}(\textbf{x}) \bigg\vert \cP_{\lam} \right] \right]\nn\\
& = &    E^o \left[   \sum\limits_{\textbf{x} \in \cP_{\lam,\neq}^n}    E^o\left[\prod\limits_{\stackrel{1 \leq j \leq k+1}{j \mbox{ odd}}} \prod\limits_{l\in B_j^e}  \left(\f{\eta  W_{x_{l-1}}W_{x_l}}{|x_l-x_{l-1}|^\al} \wedge 1\right)1_{A_j}(\textbf{x}) 
 \prod\limits_{\stackrel{1 \leq j \leq k+1}{j \mbox{ even}}} \prod\limits_{l\in B_j^e}  \left( \f{\eta  W_{x_{l-1}}W_{x_l}}{|x_l-x_{l-1}|^\al} \wedge 1\right) 1_{A_j}(\textbf{x}) \bigg \vert \cP_{\lam} \right] \right].\nn\\
\lab{eq:thetaTilde_2}
\end{eqnarray}
By applying the Cauchy-Schwarz inequality followed by using the independence of the weights in the alternating blocks, the conditional expectation on the right hand side of (\ref{eq:thetaTilde_2}) can be bounded from above by
%
%\begin{eqnarray}
%
%\lefteqn{ E^o\left[\prod\limits_{\stackrel{1 \leq j \leq k+1}{j \mbox{ odd}}} \prod\limits_{l\in B_j}  \left(\f{\eta  W_{x_{l-1}}W_{x_l}}{|x_l-x_{l-1}|^\al} \wedge 1\right) \prod\limits_{\stackrel{1 \leq j \leq k+1}{j \mbox{ even}}}^{k+1} \prod\limits_{l\in B_j}  \left( \f{\eta  W_{x_{l-1}}W_{x_l}}{|x_l-x_{l-1}|^\al} \wedge 1\right) \bigg \vert \cP_{\lam} \right] }\nn\\
%
%& \leq &    
\begin{equation}
\prod\limits_{j=1}^{k+1} 1_{A_j}(\textbf{x})  \prod\limits_{l\in B_j^e}  \left( E^o \left[ \left(\f{\eta  W_{x_{l-1}}W_{x_l}}{|x_l-x_{l-1}|^\al}\right)^2 \wedge 1 \bigg \vert \cP_{\lam} \right]\right)^{\f{1}{2}}.
%
%\prod\limits_{\stackrel{1 \leq j \leq k+1}{j \mbox{ odd}}} 1_{A_j}(\textbf{x})  \prod\limits_{l\in B_j^e}  \left( E^o \left[ \left(\f{\eta  W_{x_{l-1}}W_{x_l}}{|x_l-x_{l-1}|^\al}\right] \right)^2 \wedge 1 \bigg \vert \cP_{\lam} \right)^{\f{1}{2}}
%
% \prod\limits_{\stackrel{1 \leq j \leq k+1}{j \mbox{ even}}} 1_{A_j}(\textbf{x})  \prod\limits_{l\in B_j^e}  \left( E^o\left[ \left( \f{\eta  W_{x_{l-1}}W_{x_l}}{|x_l-x_{l-1}|^\al}\right] \right)^2 \wedge 1 \bigg \vert  \cP_{\lam} \right)^{\f{1}{2}}.
%
\lab{eq:thetaTilde_3}
\end{equation}
By Lemma 4.3 of \citep{Deijfen2013} we have
\begin{equation}
 \left[E^o \left[ \left( \f{ W_1 W_2}{t}\right)^2 \wedge 1\right]\right]^{\f{1}{2}} \leq 1_{\{t<1\}}+ c\, 1_{\{t\geq 1\}}\left( 1+(\be \vee 2) \log t \right) t^{-(\f{\be}{2} \wedge 1)} =: h(t) \mbox{ (say),}
\lab{eq:Deijfen}
\end{equation}
where $ c = \left(1+1_{\{\be\neq 2\}}\f{2}{\be-2}\right)^{\f{1}{2}} $. Using the bound from (\ref{eq:Deijfen}) in (\ref{eq:thetaTilde_3}) and then substituting in (\ref{eq:thetaTilde_2}) we obtain 
\begin{equation}
E^o \left[   \sum\limits_{\textbf{x} \in \cP_{\lam,\neq}^n} 1\{B(\textbf{x}) = B\}\right] 
\leq E^o \left[  \sum\limits_{\textbf{x} \in \cP_{\lam,\neq}^n}   \prod\limits_{j=1}^{k+1} 1_{A_j}(\textbf{x})  \prod\limits_{l\in B_j^e} h(\eta ^{-1} |x_l-x_{l-1}|^\al)   \right].
\lab{eq:thetaTilde_4}
\end{equation}
By the Campbell-Mecke formula the right hand side in (\ref{eq:thetaTilde_4}) equals
\begin{equation}
\lam^n \int \int \ldots \int  \prod\limits_{j=1}^{k+1} 1_{A_j}(\textbf{x})  \prod\limits_{l\in B_j^e} h(\eta ^{-1} |x_l-x_{l-1}|^\al)  \prod\limits_{i=1}^{n} dx_i . 
%
%\lam^n \int \int \ldots \int    \prod\limits_{\stackrel{1 \leq j \leq k+1}{j \mbox{ odd}}} \prod\limits_{l\in eB_j} h(\eta ^{-1} |x_l-x_{l-1}|^\al) \prod\limits_{\stackrel{1 \leq j \leq k+1}{j \mbox{ even}}} \prod\limits_{l\in eB_j} h(\eta ^{-1} |x_l-x_{l-1}|^\al) \prod\limits_{l=1}^{n} dx_l. 
%
\lab{eq:thetaTilde_5}
\end{equation}
We evaluate the contribution to (\ref{eq:thetaTilde_5}) from blocks of various sizes using calculations similar to those in (\ref{eq:theta_tilde_4})--(\ref{eq:theta_tilde_7}). 
%in Theorem~(\ref{thm:PhaseTransition})~(\ref{np_ercm}). 
The contribution from a block of size two equals
\begin{equation}
\int_{\bR^2}  h(\eta ^{-1} |y_1-y_2|^\al)\, dy_2 = 2\pi \int_{0}^{\infty} r\, h(\eta^{-1}\, r^\al)\, dr.
\lab{eq:thetaTilde_6}
\end{equation} 
For $m\geq 1$ the contribution from the block of size $2m+2$ is (see~(\ref{eq:theta_tilde_6a})) is
\begin{equation}
 2^{m+1}\pi \left(\int_{0}^{\infty} r^2\, h(\eta^{-1}\, r^\al)\, dr\right)  \left(\int_{0}^{\infty} r^3\, h(\eta^{-1}\, r^\al)\, dr\right)^{m-1} \left(\int_{0}^{\infty} r^2\, h(\eta^{-1}\, r^\al)\, dr\right).
 \lab{eq:thetaTilde_6a}
 \end{equation}
Using the contributions from each block $B_1,B_2, \ldots, B_{k+1}$ from (\ref{eq:thetaTilde_6}), (\ref{eq:thetaTilde_6a}) the expression in (\ref{eq:thetaTilde_5}) equals
%
%\begin{eqnarray}
%\lefteqn{\int \int \ldots \int    \prod\limits_{\stackrel{1 \leq j \leq k+1}{j \mbox{ odd}}} \prod\limits_{l\in B_j} h(\eta ^{-1} |x_l-x_{l-1}|^\al) \prod\limits_{\stackrel{1 \leq j \leq k+1}{j \mbox{ even}}} \prod\limits_{l\in B_j} h(\eta ^{-1} |x_l-x_{l-1}|^\al) \prod\limits_{l=1}^{n} dx_l }\nn\\
%
%& = & 
\begin{equation}
\left( \prod\limits_{m=0}^{\lfloor \f{n}{2}\rfloor} (2^{m+1}\pi)^{k_m} \right)  \left(\int_{0}^{\infty} r\,   h(\eta ^{-1} r^\al) \, dr\right)^{k_1} \left(\int_{0}^{\infty} r^2  h(\eta ^{-1} r^\al) \, dr\right)^{2(k-k_1)} \left(\int_{0}^{\infty} r^3  h(\eta ^{-1} r^\al) \, dr\right)^{\bar{k}},
%{\sum\limits_{\stackrel{1 \leq j \leq n}{ |B_j| \geq 6 }} \left( \f{|B_j|}{2}-1 \right) }.\nn\\
\lab{eq:thetaTilde_9}
\end{equation}
where $k_m$ and $\bar{k}$ are as defined in (\ref{eq:Kp}) and (\ref{eq:k4_eRCM}). Using (\ref{eq:thetaTilde_5}), (\ref{eq:thetaTilde_9}) in (\ref{eq:thetaTilde_2}) and substituting the resulting expression in (\ref{eq:thetaTilde_1})  we obtain
\begin{equation}
\tilde{\th}^e(\lam) \leq \lam^n \sum\limits_{k=0}^{n-1} \sum\limits_{B \in \cB_k}  \left( \prod\limits_{m=0}^{\lfloor \f{n}{2}\rfloor} (2^{m+1}\pi)^{k_m} \right) C_1^{4n}  ,
\lab{eq:thetaTilde_10}
\end{equation}
where $C_1 = \max\left\{ \int_{0}^{\infty} r^j\,   h(\eta ^{-1} r^\al) \, dr, j=1,2,3\right\}$. By the same arguments as in Theorem~\ref{thm:PhaseTransition}~(\ref{np_ercm}) there exists a constant $C$ such that
\[ \tilde{\th}^e(\lam) \leq  (C \lam)^n \to 0 \]
as $n \to \infty$ for all $\lam \in (0, \lam_0)$, for some $\lam_0 > 0$ provided $ \int_{0}^{\infty} r^3  h(\eta ^{-1} r^\al) \, dr <\infty$, which from (\ref{eq:Deijfen}) is true since $\al>4$ and  $\al\be > 8$. \qed

\subsection{Proof of Theorem~\ref{thm:RSW_nonpercolation} \, for ieRCM}
\lab{section:RSW_for_ieRCM}
The RSW results for the ieRCM as enumerated in Theorem~\ref{thm:RSW_nonpercolation} follow in a manner identical to that for the eRCM once we prove the analog of Proposition~\ref{prop:largest_edge_length_in_box} for the length of the longest edge in the graph $H_{\lam}$. It thus suffices to prove the following proposition, the proof of which is identical to that of Proposition~\ref{prop:largest_edge_length_in_box} with the obvious changes.
\begin{prop}
Let $\min\{\al,\al\be\} > 4$ and consider the graph $H_\lam$ with connection function $g$ satisfying defined as in (\ref{eq:connection_function}). For any $s > 0$ let $M_{s}$ be the length of the longest edge in $H_{\lam}$ intersecting the box $B_{s} = [-s,s]^2$. Then for any $t > 0$ and $\tau>\f{2}{\min\{\al,\al\be\}-2}$ we have $P\left(M_{ts}>s^{\tau}\right)\rar 0$ as $ s\rar \infty$.
\lab{prop:largest_edge_length_IN}
\end{prop}

We shall use the following upper bound on the expected value of the connection function. The proof will be given later. We shall write $E^{W_x, W_y}$ to denote the expectation with respect to the random weight $W_x$ and $W_y$.
\begin{lem} For $\al, \be, \eta >0$ and any $x\in \bR^2$ with $|x|^{\al} > \eta$ there exists constants $C_1, C_2, C_3$ such that the connection function given by (\ref{eq:connection_function}) satisfies
\begin{equation}
    E^{W_o,W_x}[g(O,x)]  = 
\begin{cases}
    C_1|x|^{-\al}+C_2|x|^{-\al\be}\log|x|+C_3|x|^{-\al\be},                     & \text{if } \be \neq 1\\
    C_1|x|^{-\al}(\log|x|)^2+C_2|x|^{-\al}\log|x|+C_3|x|^{-\al},              & \text{if } \be =1.
\end{cases}
\lab{eq:upper_bound_g} 
\end{equation}
\lab{lemma:upper_bound_g}
\end{lem} 
{\bf Proof of Proposition~\ref{prop:largest_edge_length_IN}.}   Fix $c > 4$, $t>0$. Let $B(O,s):=\{x\in \bR^2:|x|\leq s\}$ be the ball of radius $s$ centered at the origin. Define the events $D_s(l)=\{M_s>l \}$,
$$O_{s}(\tau)=\{ X \in \cP_{\lam}: \mbox{ there is an edge of length longer than } s^{\tau} \mbox{ incident on } X \mbox{ in } H_{\lam} \},$$ 
$$\bar{O}_{t,s}(\tau)=\{ (X, Y) \in \cP_{\lam}^2: \mbox{ there is an edge in } H_{\lam} \mbox{ joining } X,Y \mbox{ of length longer than } s^{\tau}, \overline{X Y} \mbox{ intersects } B(O, \sqrt{2}ts)\}.$$ 

Recall that for any two points $x,y \in \bR^2$, $\overline{xy}$ denotes the line segment joining $x$ and $y$. 
\begin{eqnarray}
P\left(D_{ts}(s^{\tau})\right)&\leq & E\left[\sum_{X,Y\in \cP_{\lam}}1\left\{\overline{X Y}\mbox{ intersects } B_{ts}\right\} 1\left\{|X- Y|\geq s^\tau\right\}\right]\lab{eq:1st_inequality_i}\nn\\
& \leq &  E\left[\sum_{X\in \cP_{\lam}\cap B(O,\sqrt{2}ts)}1\left\{X\in O_{s}(\tau)\right\}\right] + 
E \left[\sum_{X,Y \in \cP_{\lam}\cap B(O,\sqrt{2}ts)^c }1\left\{(X,Y) \in \bar{O}_{t,s}\right\}\right].
\lab{eq:2nd_inequality_IN}
\end{eqnarray}
In what follows we shall write $E^{W_x}$ to denote the expectation with respect to the random weight $W_x$. By the Campbell-Mecke formula applied to the first term on the right hand side of the last inequality in (\ref{eq:2nd_inequality_IN}) we obtain
\begin{eqnarray}
E \left[\sum_{X\in \cP_{\lam}\cap B(O,\sqrt{2}ts)}1\left\{X\in O_{s}(\tau)\right\}\right]
&=& C \lam (ts)^{2}E^{W_o}\left[ P^o\left(O\in O_{s}(\tau)\vert W_o\right) \right] \nn\\
& = & C \lam (ts)^{2}E^{W_0}\left[ 1-P^o\left(\mbox{none of the edges incident on O is of length} \geq s^{\tau} \vert W_o \right)\right]\nn\\
& = & C \lam (ts)^2 E^{W_0} \left[ 1- \exp{\left(-\lam \int_{B(O, s^{\tau})^c}E^{W_x}\left[g(|x|) \big\vert W_o \right] dx\right)}\right] \nn\\
 & \leq & C (\lam ts)^2 E^{W_0} \left[ \int_{B(O, s^{\tau})^c}E^{W_x}\left[g(|x|) \big\vert W_o \right] dx \right] \nn\\
  & \leq & C (\lam ts)^2  \int_{B(O, s^{\tau})^c}E^{W_o,W_x}\left[1-\exp\left( -\f{\eta W_oW_x}{|x|^\al}\right) \right] dx ,
\lab{eq:one_in_IN}
\end{eqnarray}
where we have used the fact that conditional on the weight $W_o$ at the origin $O$, the points of $\cP_{\lam}$ from which there is incident on $O$ an edge that is of length longer than $s^{\tau}$ is a Poisson point process of intensity $\lam E^{W_x}\left[g(|x|) \vert W_o \right] 1\{x \in B(O, s^{\tau})^c\}$ and the inequality $1-e^{-y}\leq y$. Using Lemma~\ref{lemma:upper_bound_g} and the fact that  $\log r \leq r^\ep$,
\begin{eqnarray}
\int_{B(O, s^{\tau})^c}E^{W_o,W_x}\left[1-\exp\left( -\f{\eta W_oW_x}{|x|^\al}\right) \right] dx   & = & C_1 s^2 \int_{s^\tau}^{\infty}r \left( C_1 r^{-\al}+C_2 r^{-\al\be}\log r+C_3 r^{-\al\be}\right)\,dr \nn\\
%
%& \leq & C_1 s^2 \int_{s^\tau}^{\infty}r \left( C_1 r^{-\al}+C_2 r^{-\al\be} r^\ep +C_3 r^{-\al\be}\right)\,dr \nn\\
%
& \leq & C_4\, s^{2-\tau(\al-2)}+C_5\, s^{2-\tau(\al\be-2)}+C_6\, s^{2-\tau(\al\be-\ep -2)},
 \lab{eq:}
 \end{eqnarray}
  Similarly we can bound the second term on the right hand side in the last inequality in (\ref{eq:2nd_inequality_IN}) as follows (see Figure~\ref{fig:both_out}). Let $\tilde{g}(x, y) := E[g(x,y)]$, $x, y\in \bR^2$, where $g$ is as specified in (\ref{eq:connection_function}). Observe that $\tilde{g}(x,y)$ depends on $x,y$ only via $|x-y|$ and so by an abuse of notation we will write $\tilde{g}(|x-y|)$ for $\tilde{g}(x,y)$.
\begin{eqnarray}
 E \left[\sum_{X,Y \in \cP_{\lam}\cap B(O,\sqrt{2}ts)^c }1\{(X,Y) \in \bar{O}_{t,s}\}\right] & \leq &  \lam^2  \int_{ B(O,\sqrt{2}ts)^c} \int_{ D_x \cap B(x, s^{\tau})^c} E^{W_x,W_y}[g(x,y)]\,dy\,dx \nn\\
 & = &  \lam^2  \int_{ B(O,\sqrt{2}ts)^c} \int_{ D_x \cap B(x, s^{\tau})^c} \tilde{g}(|x-y|)\,dy\,dx,
\lab{eq:both_out_1_IN}
\end{eqnarray}
where $D_x$ is the unbounded region $ATS'T'A'$ as shown in Figure~\ref{fig:both_out}. Changing to polar coordinates as in the proof of Proposition~\ref{prop:largest_edge_length_in_box} and using Lemma~\ref{lemma:upper_bound_g} we can bound the last  expression in (\ref{eq:both_out_1_IN}) by
\begin{eqnarray}
\lefteqn{C_7 \int_{\sqrt{2}ts}^{\infty}\int_{ 0}^{2\pi}R\,d\phi\,dR\int_{\left(s^\tau \vee \sqrt{R^2-2t^2s^2}\right)}^{\infty}r \, \tilde{g}(r)\,dr}\nn\\
 &\leq& C_8 \int_{\sqrt{2}ts}^{\infty}\int_{\left(s^\tau \vee \sqrt{R^2-2t^2s^2}\right)} ^{\infty} r\left(  r^{-\al}+ r^{-\al\be}\log r+ r^{-\al\be}\right)\, dr R\,dR\nn\\
&\leq& C_9 \int_{\sqrt{2}ts}^{\infty}  \left( \left(s^\tau \vee \sqrt{R^2-2t^2s^2}\right)^{2-\al}+ \left(s^\tau \vee \sqrt{R^2-2t^2s^2}\right)^{2-\al\be}\right) R\,dR.
\lab{eq:both_out_2_i}
\end{eqnarray}
%
%The first step follows from the Lemma~\ref{lemma:upper_bound_g}, $\tilde{g}(r)\leq C_1 r^{-\al}+C_2 r^{-\al\be}\log r+C_3 r^{-\al\be}$. 
The integral of the first integrand on right hand side in (\ref{eq:both_out_2_i}) can be evaluated as follows.
\begin{eqnarray}
 \int_{ \sqrt{2} ts}^{\infty} \left(s^\tau \vee \sqrt{R^2-2t^2s^2}\right)^{2-\al}R\,dR
 &=&\int_{\sqrt{2}ts}^{\sqrt{2t^2s^2+s^{2\tau}}}s^{\tau(2-\al)}R\,dR+ \int_{\sqrt{2t^2s^2+s^{2\tau}}}^{\infty} \left(R^2-2t^2s^2\right)^{\f{2-\al}{2}}R\,dR\nn\\
&=& C_{10} s^{-\tau(\al-4)}.
%C_2 s^{-\tau(\al-4)} +C_3\int_{s^\tau}^{\infty} u^{3-\al}\,du=C_4s^{-\tau(\al-4)}
\lab{eq:both_outside_4}
\end{eqnarray}
Similarly the second term on the right hand side in (\ref{eq:both_out_2_i}) can be evaluated to obtain 
\begin{eqnarray}
\int_{ \sqrt{2}ts}^{\infty} \left(s^\tau \vee \sqrt{R^2-2t^2s^2}\right)^{2-\al\be}R\,dR &=& C_{11} s^{-\tau(\al\be-4)}.
\lab{eq:both_outside_5}
\end{eqnarray}
Substituting  (\ref{eq:both_outside_4}), (\ref{eq:both_outside_5}) in (\ref{eq:both_out_2_i}) and using (\ref{eq:one_in_IN}) and (\ref{eq:both_out_2_i})  in (\ref{eq:1st_inequality_i}) we obtain
\begin{equation}
P\left(D_{ts}(s^{\tau})\right)\leq C_{12} \left( s^{2-\tau(\al-2)} + s^{2-\tau(\al\be-2)} + s^{-\tau(\al-4)} + s^{-\tau(\al\be-4)} \right) \rar 0, \nn
\end{equation}
as $s\rar \infty$ since $\tau>\f{2}{\min\{\al,\al\be\}-2}$ and $\min\{\al,\al\be\} > 4$. \qed

\textbf{Proof of Lemma~\ref{lemma:upper_bound_g}.} We will prove the result for the case $\be \ne 1$. The proof for $\be = 1$ follows with minor changes. Fix $x\in \bR^2$ such that $|x|^{\al} > \eta$. Then
\begin{align}
E^{W_o,W_x}[g(O,x)] &= E^{W_o,W_x}\left[1-\exp\left(-\f{\eta W_0W_x}{|x|^\al}\right)\right]\nn\\
&\leq E^{W_o,W_x}\left[\f{\eta W_0W_x}{|x|^\al}\wedge 1\right]\nn\\
&=E^{W_o,W_x}\left[\f{\eta W_0W_x}{|x|^\al}; W_0W_x<\f{|x|^\al}{\eta}\right]+P\left(W_0W_x\geq \f{|x|^\al}{\eta}\right).
\lab{eq:two_terms}
\end{align}
By our assumption on the distribution of the weights it is easy to see that the product $W_0W_x$ has a density given by
\begin{equation}
f(w)=\be^2w^{-\be-1}\log w,\,  w\geq 1.
\lab{eq:density_w0wx}
\end{equation} 
The first term on the right in (\ref{eq:two_terms}) can be evaluated using the expression in (\ref{eq:density_w0wx}) and the integration by parts formula to yield
\begin{eqnarray}
E^{W_o,W_x}\left[\f{\eta W_0W_x}{|x|^\al}; W_0W_x<\f{|x|^\al}{\eta}\right] & = &  \f{\eta}{|x|^\al} \int_{1}^{|x|^\al / \eta}\, w\,f(w)\, dw \nn \\
 & \leq &  c_1|x|^{-\al}+c_2|x|^{-\al\be}\log|x|+c_3|x|^{-\al\be},
\lab{eq:first_term_1}
\end{eqnarray}
for some constants $ c_1,c_2,c_3$. The second term on the right in (\ref{eq:two_terms}) can be computed in a similar fashion.
\begin{eqnarray}
P\left(W_0W_x\geq \f{|x|^\al}{\eta}\right) & = &  \int_{|x|^\al / \eta}^{\infty}f(w)\, dw \nn\\
& \leq &   c_4|x|^{-\al\be}+c_5|x|^{-\al\be}\log |x|,
\lab{eq:second_term_1}
\end{eqnarray}
for some constants $c_4,c_5$. (\ref{eq:upper_bound_g}) now follows by substituting from (\ref{eq:first_term_1}),(\ref{eq:second_term_1}) in (\ref{eq:two_terms}). \qed
\subsection{Proof of Theorem~\ref{thm:PhaseTransition}~(\ref{np_ps})}
%
%Let $L_o$ be the stick with mid point at the origin $(O)$ and be of length $2\ell$. Let $R$ be the half length of a stick $L_1$ that intersects $L_o$ at an angle $\th$. The following proposition provides the size-biased density corresponding to the random variable $R$ given $\th$ and $\ell$. 

To show that $\lam_{PS} < \infty$ consider the graph $PS_\lam$ with half length density $h$ satisfying $0<\int_{0}^{\infty}\ell^2 \,h(\ell)\,d\ell<\infty$.  Let $R_0:=\inf\{\ell : h(\ell) > 0 \}$. Pick any $R_1 > R_0$ finite such that $\int_{R_0}^{R_1} h(\ell)\,d\ell > 0$.  Consider the graph $PS_\lam$ and set sticks of length greater than $2R_1$ to be equal to $2R_1$ to obtain the graph $\tilde{PS}_\lam$. We have from \cite{RRoy1991} that the critical threshold parameter $\tilde{\lam}_{PS} < \infty$ for the truncated model. Since $PS_\lam$ percolates if $\tilde{PS}_\lam$ does, we have $\lam_{PS} <\infty$.

\begin{figure}[h!]
 \centering
\includegraphics[width=0.7\linewidth]{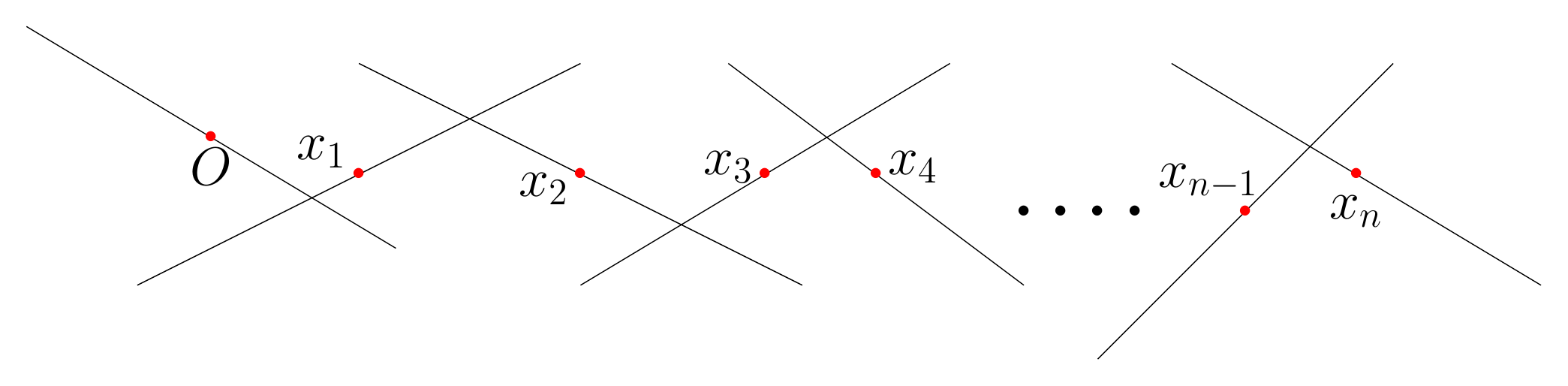}
  \captionof{figure}{A self-avoiding path in $PS_{\lam}$}
  \label{fig:PS_Path}
\end{figure}

The proof of $\tilde{\lam}_{PS}> 0$ is simpler for the Poisson stick model since for each ordered sequence $\textbf{x} = (x_1, x_2, \cdots, x_n)$, there is only way in which a path can occur (see Figure~\ref{fig:PS_Path}). For $x \in \cP_{\lam}$ let $L(x)$ denote the stick centered at $x$. We will denote $L(x_j)$ by $L_j$. Then
\begin{eqnarray}
\th_{PS}(\lam) &\leq &  P^o(\mbox{ there is a self-avoiding  path on } n \mbox{ vertices in } PS_\lam \mbox{ starting from } O)\nn\\
& \leq &  E^o \left[\sum\limits_{\textbf{x} \in \cP_{\lam,\neq}^n} 1\{O \to x_1 \to x_2 \to \cdots \to x_n \mbox{ occurs }\}     \right] \nn\\
& = &     E^o \left[   \sum\limits_{\textbf{x} \in \cP_{\lam,\neq}^n} \prod\limits_{j=1}^{n} 1\{L_j \mbox{ intersects } L_{j-1}\}   \right].
\lab{eq:thetaPS_1}
\end{eqnarray}
\begin{figure}[h!]
 \centering
\includegraphics[width=0.6\linewidth]{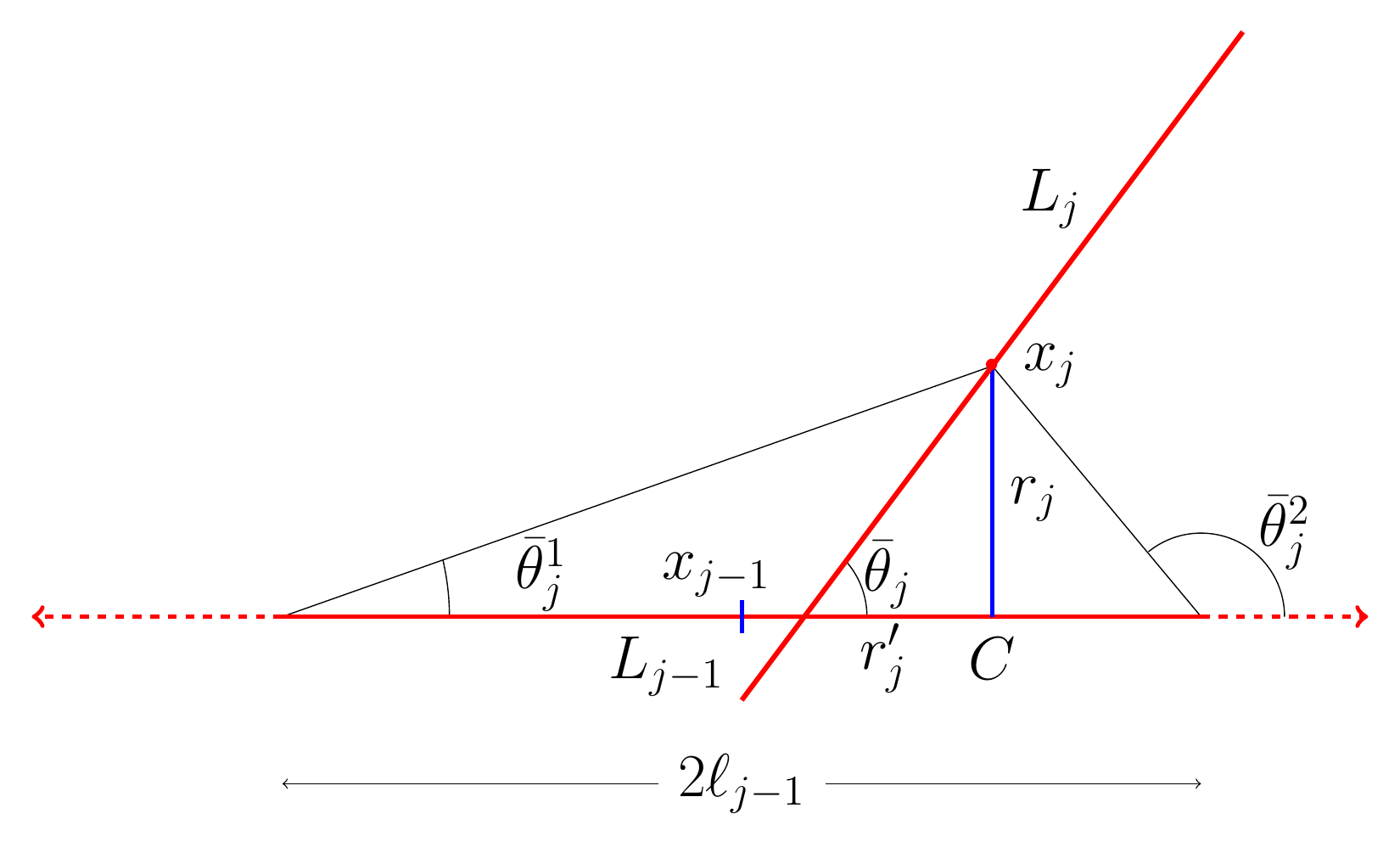}
  \captionof{figure}{The stick $L_j$ intersects $L_{j-1}$}
  \label{fig:size_biased_measure}
\end{figure}
Let $\th_j$ be the orientation of the stick $L_j$, $\bar{\th}_j$ be the relative orientation of $L_j$ with respect to $L_{j-1}$ and $F_{\bar{\th}_{j-1}}(\bar{\th}_j)$ be the distribution of $\bar{\th}_j$ given $\th_{j-1}$. Let $(r_j',r_j)$ be the coordinates of $x_j$ with $L_{j-1}$ being the horizontal axis and $x_{j-1}$ as the origin. Given $x_{j-1}, x_j, \ell_{j-1}$, $L_j$ intersects $L_{j-1}$ provided $\bar{\th}_j \in \left[\bar{\th}_j^1 \wedge \bar{\th}_j^2, \bar{\th}_j^1 \vee \bar{\th}_j^2 \right]$ and $\ell_j \geq \f{r_{j}}{\sin\bar{\th}_j}$ (see Figure~\ref{fig:size_biased_measure}). We apply the Campbell-Mecke formula and then change the variable $x_j$ to $(r_j',r_j)$. The expectation in (\ref{eq:thetaPS_1}) then evaluates to 
\begin{equation} 
\lam^n \int_{0}^{\pi} \int_{0}^{\infty}   h(l_0)\, dF(\th_0)\, d\ell_0 \, \prod_{j=1}^n \int_{-\infty}^{\infty} dr_j \int_{-\infty}^{\infty} dr'_j \int_{\bar{\th}_j^1 \wedge \bar{\th}_j^2}^{\bar{\th}_j^1 \vee \bar{\th}_j^2} dF_{\bar{\th}_{j-1}}(\bar{\th}_j) \int_{\f{r_{j}}{\sin\bar{\th}_j}}^{\infty} h(\ell_j) \, d\ell_j. 
\lab{eq:thetaPS_1a}
\end{equation} 
By Fubini's Theorem, exchanging the integrals with respect to $r'_j$ and $\bar{\th}_j$ the expression in 
(\ref{eq:thetaPS_1a}) is bounded by
\begin{eqnarray}
 \lefteqn{\lam^n \int_{0}^{\pi} \int_{0}^{\infty}   h(l_0)\, dF(\th_0) \, d\ell_0\, \prod_{j=1}^n \int_{-\infty}^{\infty} dr_j \int_0^{\pi} dF_{\bar{\th}_{j-1}}(\bar{\th}_j)   \int_{-\ell_{j-1}+ r_j\tan\bar{\th}_j}^{\ell_{j-1} + r_j\tan\bar{\th}_j} dr'_j \int_{\f{r_{j}}{\sin\bar{\th}_j}}^{\infty} h(\ell_j) \, d\ell_j } \nn\\
 &=&\lam^n \int_{0}^{\pi} \int_{0}^{\infty}   h(l_0)\, dF(\th_0) \, d\ell_0 \, \prod_{j=1}^n 2\ell_{j-1} \int_{-\infty}^{\infty} dr_j \int_0^{\pi} dF_{\bar{\th}_{j-1}}(\bar{\th}_j)    \int_{\f{r_{j}}{\sin\bar{\th}_j}}^{\infty} h(\ell_j) \, d\ell_j
\lab{eq:thetaPS_1b}
\end{eqnarray}
Consider the integrals with respect to the index $j=n$ in the expression on the right in (\ref{eq:thetaPS_1b}). By another use of Fubini's Theorem we obtain 
\begin{eqnarray}
\int_{-\infty}^{\infty} dr_n \int_0^{\pi} dF_{\bar{\th}_{n-1}}(\bar{\th}_n)    \int_{\f{r_{n}}{\sin\bar{\th}_n}}^{\infty} h(\ell_n) \, d\ell_n & = &  \int_{0}^{\pi} dF_{\bar{\th}_{n-1}}(\bar{\th}_n) \int_{0}^{\infty} h(\ell_n) \, d\ell_n   \int_{0}^{\ell_n \sin\bar{\th}_n} dr_n \nn \\
 & = & \int_{0}^{\pi} \sin\bar{\th}_n \; dF_{\bar{\th}_{n-1}}(\bar{\th}_n) \int_{0}^{\infty}\ell_n  h(\ell_n) \, d\ell_n  \leq \int_{0}^{\infty}\ell  h(\ell) \, d\ell.
\lab{eq:thetaPS_2}
\end{eqnarray}
Substituting from (\ref{eq:thetaPS_2}) in (\ref{eq:thetaPS_1b}) yields the bound
\begin{equation}
\lam^n \int_{0}^{\pi} \int_{0}^{\infty}   h(l_0)\, dF(\th_0) \, d\ell_0 \prod_{j=1}^{n-1} 2\ell_{j-1} I_{j}   \left( \int_{0}^{\infty}\ell  h(\ell) \, d\ell  \right), 
\lab{eq:thetaPS_3}
\end{equation}
where 
\begin{eqnarray}
I_j &:= &\int_{0}^{\pi} dF_{\bar{\th}_{j-1}}(\bar{\th}_j) \int_{0}^{\infty}dr_j \int_{\f{r_j}{\sin\bar{\th}_j}}^{\infty} \ell_j \, h(\ell_j) \, d\ell_j \nn\\
&=& \int_{0}^{\pi} dF_{\bar{\th}_{j-1}}(\bar{\th}_j)  \int_{0}^{\infty}\ell^2_j \sin\bar{\th}_j h(\ell_j) \, d\ell_j  
 \leq  \int_{0}^{\infty}\ell^2  h(\ell) \, d\ell.
\lab{eq:thetaPS_4}
\end{eqnarray}
%
%Substituting from (\ref{eq:thetaPS_4})  in (\ref{eq:thetaPS_3})  we obtain an upper bound as 
%
%\begin{equation}
%
%  \lam^n \int_{0}^{\pi} \int_{0}^{\infty}   h(l_0)\, dF(\th_0)  \prod_{j=1}^{n-2} 2\ell_{j-1} I_{j} 2^2 \left( \int_{0}^{\infty}\ell^2  h(\ell) \, d\ell  \right) \left( \int_{0}^{\infty}\ell  h(\ell) \, d\ell  \right).\nn
%
%\end{equation}
%
By iteratively evaluating the $I_j$ and using (\ref{eq:thetaPS_4}) we obtain
\begin{eqnarray*} 
\th_{PS}(\lam) & \leq & 2^{n-1}\lam^n \int_{0}^{\pi} \int_{0}^{\infty}   h(l_0)\, dF(\th_0) \, d\ell_0 \, 2 \ell_0 \left(\int_{0}^{\infty} \ell^2\,  h(\ell) \, d\ell \right)^{n-2}  \int_{0}^{\infty} \ell\,  h(\ell) \, d\ell \nn\\
&=&  2^n \lam^n   \left(\int_{0}^{\infty} \ell\,  h(\ell) \, d\ell\right)^2 \left(\int_{0}^{\infty} \ell^2\,  h(\ell) \, d\ell\right)^{n-2} \to 0
\end{eqnarray*}
as $n \to \infty$ for all $\lam  > 0$ sufficiently small provided $ \int_{0}^{\infty} \ell^2\,  h(\ell) \, d\ell < \infty$. This completes the proof of the Theorem~\ref{thm:PhaseTransition}~(\ref{np_ps}). \qed

\subsection{Proof of Theorem~\ref{thm:RSW_nonpercolation} for Poisson stick model }
\lab{section:RSW_for_PS}
The proof follows along the same lines as that for the eRCM using the following Proposition on the length of the longest stick in $PS_{\lam}$. Since $c > 3$ we can choose a $\tau < 1$ so that the conclusion of Proposition~\ref{prop:largest_stick_length} holds, which gives us the condition under which the rest of the proof works. 

\begin{figure}
\centering
\includegraphics[width=0.7\linewidth]{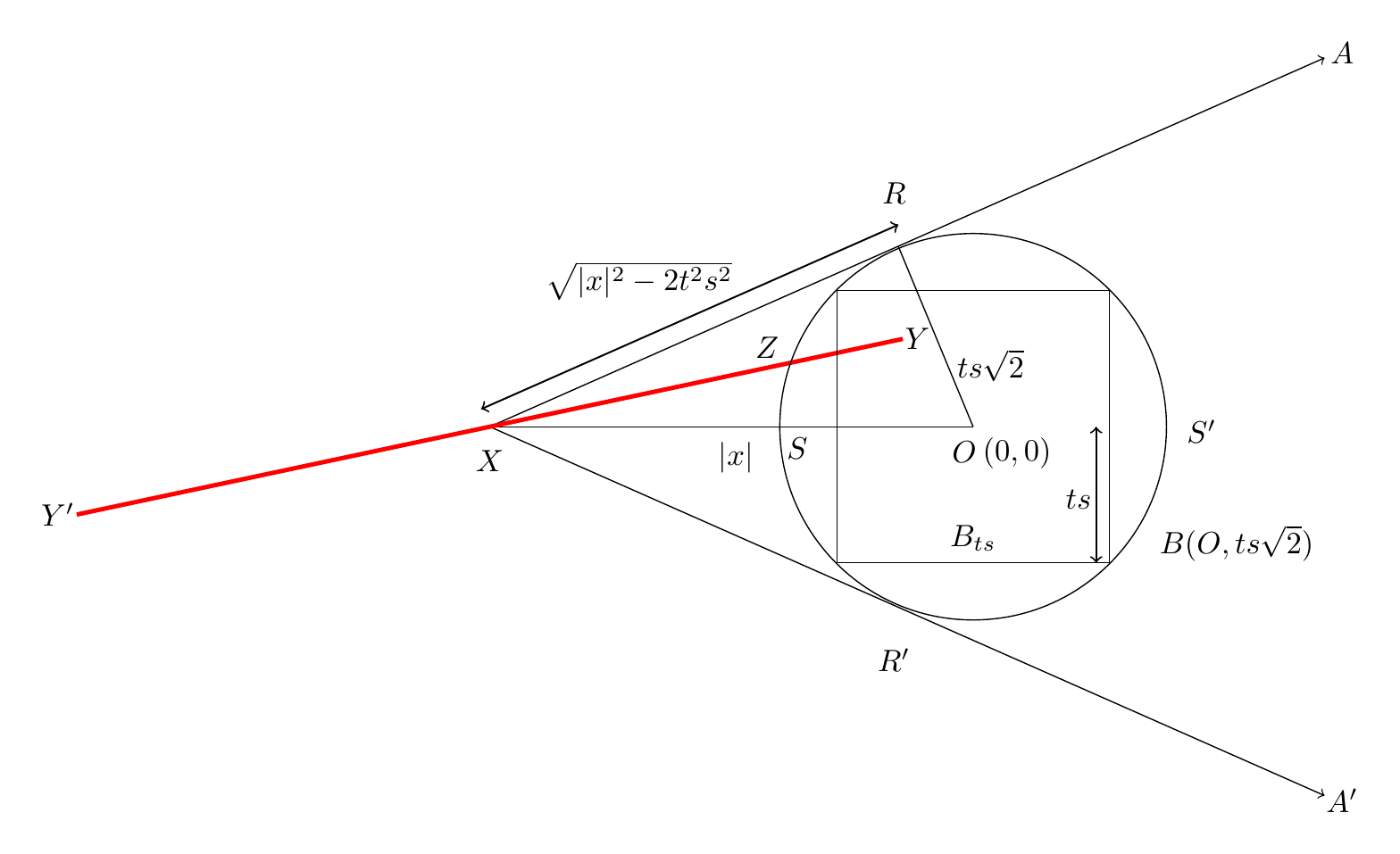}
  \captionof{figure}{The red stick with its mid point outside the ball $B(O, \sqrt{2}ts)$.}
  \label{fig:out_stick_intersects}
\end{figure}

\begin{prop}
Consider the graph $PS_\lam$ with stick's half length density function $h$ satisfying $h(\ell) = O(\ell^{-c})$ as $\ell \rar \infty$ for some $c>3$. For any $s > 0$ let $\bar{M}_{s}$ be the half length of the longest stick in $PS_{\lam}$ intersecting the box $B_{s} = [-s,s]^2$. Then for any $t > 0$ and $\tau>\f{2}{c-1}$ we have $P\left(\bar{M}_{ts}>s^{\tau}\right)\rar 0$ as $s\rar \infty$.
\lab{prop:largest_stick_length}
\end{prop}
{\bf Proof of Proposition~\ref{prop:largest_stick_length}.} Fix $c > 3$, $t>0$ and $\tau>\f{2}{c-1}$. Since $c > 3$ it should suffice to let $\tau > \f{2}{c-1}$. Recall that $B(O,s)$ denotes the ball of radius $s$ centered at the origin and for $X \in \cP_{\lam}$, $L_X$ denotes the stick with mid point at $X$ with half length random variable distributed independently according to the probability density function $h$. Define the events $\bar{D}_s(\ell)=\{\bar{M}_s>\ell \}$, 
\[ O_{s}(\tau)=\{ X \in \cP_{\lam} : L_X \mbox{ has half length longer than } s^{\tau} \mbox{ and } \mbox{ intersects } B(O, \sqrt{2}ts)\}. \]
\begin{eqnarray}
P\left(\bar{D}_{ts}(s^{\tau})\right) &\leq & E\left[\sum_{X\in \cP_{\lam}}1\left\{L_X\mbox{ intersects } B_{ts}\right\}1\left\{L_X\geq s^\tau\right\}\right]\nn\\
&\leq & E\left[\sum_{X\in \cP_{\lam}\cap B(O, \sqrt{2}ts)}1\left\{X\in O_{s}(\tau)\right\}\right] +  E \left[\sum_{X\in \cP_{\lam}\cap B(O, \sqrt{2} ts)^c }1\left\{X\in O_{t,s}(\tau)\right\}\right].
\lab{eq:2nd_ineq_PS}
\end{eqnarray}
By our assumption $h(\ell)\leq C\,\ell^{-c}$ for all $\ell$ sufficiently large. The Campbell-Mecke formula applied to the first term on the right hand side of the last inequality in (\ref{eq:2nd_ineq_PS}) yields
\begin{eqnarray}
\lefteqn{E\left[\sum_{X\in \cP_{\lam}\cap B(O, \sqrt{2} ts)}1\left\{X\in O_{s}(\tau)\right\}\right] 
 =  C_0\lam (ts)^{2}P^o \left(L_O\geq s^\tau\right)}\nn\\
& = & C_1 s^2 \int_{s^\tau}^{\infty} \int_{0}^{\pi} h(\ell)\, dF(\th)\,d\ell
=C_1 s^2 \int_{s^\tau}^{\infty}h(\ell)\,d\ell \leq C_2\, s^{2-\tau(c-1)},
\lab{eq:one_in_PS}
\end{eqnarray}
for all $s$ sufficiently large. For the second term in (\ref{eq:2nd_ineq_PS}) using the assumption on $h$, the Campbell-Mecke formula and the fact that $\left|\overline{XY}\right| \geq s^\tau \vee (R-\sqrt{2} ts)$ (see Figure~\ref{fig:out_stick_intersects}) we obtain for all $s$ sufficiently large
\begin{eqnarray}
E \left[\sum_{X\in \cP_{\lam}\cap B(O,\sqrt{2} ts)^c }1\left\{X\in O_{t,s}(\tau)\right\}\right] 
& \leq &   C \lam\int_{ \sqrt{2} ts}^{\infty} R\,dR\int_{s^\tau \vee (R-\sqrt{2}ts)}^{\infty} \int_{0}^{\pi} h(\ell)\, dF(\th)\,d\ell\nn\\
& = & C_1 \int_{\sqrt{2}ts}^{\infty}\left(s^\tau \vee \left(R-\sqrt{2}ts\right)\right)^{1-c}R\,dR \nn\\
&=& C_2 \int_{ \sqrt{2}ts}^{\sqrt{2}ts + s^{\tau}} s^{\tau(1-c)}R\,dR  +  C_3 \int_{\sqrt{2}ts + s^{\tau}}^{\infty} \left(R-\sqrt{2}ts\right)^{1-c}R\,dR\nn\\
&=& C_4 s^{-\tau(c-3)} + C_5s^{-(\tau(c-2)-1)} +C_6 \int_{s^\tau}^{\infty} u^{1-c}\left(u+\sqrt{2}ts\right)\,du \nn\\
& = & C_7s^{-\tau(c-3)}+C_8s^{-(\tau(c-2)-1)}.
\lab{eq:out_PS1}
\end{eqnarray}
Substituting from (\ref{eq:one_in_PS}) and (\ref{eq:out_PS1}) in (\ref{eq:2nd_ineq_PS}) we obtain
\begin{eqnarray} 
P\left(\bar{D}_{ts}(s^\tau)\right)&\leq&  c_1 \, s^{-(\tau(c-1)-2)} + c_2 s^{-\tau(c-3)}+c_3s^{-(\tau(c-2)-1)} \to 0,\nn
\end{eqnarray}
as $s\rar \infty$, since $\tau> \f{2}{c-1} > \f{1}{c-2}$ for $c>3$.\qed
\section*{Acknowledgements:}
The authors are thankful to Mathew D. Penrose and Yogeshwaran D.  for many useful discussions and references.
%

%\bibliographystyle{plain}
%\bibliography{percolation_eRCM}

\begin{thebibliography}{10}

\bibitem{Ahlberg2018}
D.~{Ahlberg}, V.~{Tassion}, and A.~{Teixeira}.
\newblock Sharpness of the phase transition for continuum percolation in
  $\mathbb{R}^2$.
\newblock {\em Probab. Th. Rel. Fields}, 172(1):525--581, 2018.

\bibitem{Aldous2014}
D.~J. Aldous.
\newblock Scale-invariant random spatial networks.
\newblock {\em Electron. J. Probab.}, 19(15):1--41, 2014.

\bibitem{AldousShun2010}
D.~J. Aldous and J.~Shun.
\newblock Connected spatial networks over random points and a route-length
  statistic.
\newblock {\em Statist. Sci.}, 25(3):275--288, 2010.

\bibitem{AldousKendall2008}
D.~J.~{Aldous} and W.~S.~{Kendall}. 
\newblock Short-length routes in low-cost networks via Poisson line patterns.
\newblock {\em Adv. Appl. Prob.}, 40(1):1--21, 2008.

\bibitem{Alexander1996}
K. S.~Alexander.
\newblock The RSW theorem for continuum percolation and the CLT for Euclidean minimal spanning trees.
\newblock {\em Ann. Appl. Prob.},  6(2):466--494, 1996.

\bibitem{Balberg1983}
I.~Balberg and N.~Binenbaum.
\newblock Computer study of the percolation threshold in a two-dimensional
  anisotropic system of conducting sticks.
\newblock {\em Phys. Rev. B}, 28(7):3799--3812, 1983.

\bibitem{Bollobas2001}
B.~Bollob\`{a}s.
\newblock {\em Random graphs}, volume~73 of {\em Cambridge Studies in Advanced
  Mathematics}.
\newblock Cambridge University Press, Cambridge, 2001.

\bibitem{Daniels2016}
C.~J.~E. Daniels.
\newblock {On the phase transition in certain percolation models}.
\newblock{PhD Thesis, University of Bath}, 2016.

\bibitem{Deijfen2013}
M. {Deijfen}, R. v. d. {Hofstad}, and G. {Hooghiemstra}.
\newblock Scale-free percolation.
\newblock {\em Ann. Inst. H. Poincaré Probab. Statist.}, 49(3):817--838, 2013.

\bibitem{Deprez2015}
P. {Deprez}, R. S. {Hazra}, and M. V. {W{\"u}thrich}.
\newblock Inhomogeneous long-range percolation for real-life network modeling.
\newblock {\em Risks}, 3(1):1--23, 2015.

\bibitem{Deprez2018}
P. {Deprez} and M. V. {W{\"u}thrich}.
\newblock Scale-free percolation in continuum space.
\newblock {\em  Commun. Math. Stat.}, 7, 269--308, 2019.


\bibitem{Duminil-Copin2016}
H. Duminil-Copin and V. Tassion.
\newblock A new proof of the sharpness of the phase transition for Bernoulli
  percolation and the ising model.
\newblock {\em Commun. Math. Phy.}, 343(2):725--745, Apr
  2016.

\bibitem{Erdos59}
P.~Erd\"{o}s and A.~R\'{e}nyi.
\newblock On random graphs.
\newblock {\em Publ. Math. Debrecen}, 6:290--297, 1959.

\bibitem{Erdos60}
P.~Erd\"{o}s and A.~R\'{e}nyi.
\newblock On the evolution of random graphs.
\newblock {\em Publications of the Mathematical Institute of the Hungarian
  Academy of Science}, 5:17--61, 1960.

\bibitem{Franceschetti2007}
M.~Franceschetti and R.~Meester.
\newblock {\em Random Networks for Communication}.
\newblock Cambridge University Press, 2007.

\bibitem{Gilbert59}
E.~N. Gilbert.
\newblock Random Graphs.
\newblock {\em Ann. Math. Statist.}, 30(4):1141--1144, 1959.

\bibitem{Gilbert61}
E.~N. Gilbert.
\newblock Random Plane Networks.
\newblock {\em J. Soc. Indust. Appl. Math.}, 9(4):533--543, 1961.

\bibitem{Grimmett1999}
G.~Grimmett.
\newblock {\em Percolation}.
\newblock Springer, 1999.

\bibitem{Haenggi2012}
M.~Haenggi.
\newblock {\em Stochastic Geometry for Wireless Networks}.
\newblock Cambridge University Press, Cambridge, 2012.

\bibitem{Hu2004}
L.~{Hu} and G.~{Hecht}, D. S.~and{ Gr\"{u}ner}.
\newblock Percolation in transparent and conducting carbon nanotube networks.
\newblock {\em Nano Lett.}, 4(12):2513--2517, 2004.

\bibitem{Janson2000}
S.~{Janson}, T.~{\L uczak}, and A.~{Rucinski}.
\newblock {\em Random Graphs}.
\newblock Wiley-Interscience Series in Discrete Mathematics and Optimization.
  2000 John Wiley \& Sons, 2000.

\bibitem{Kesten1980}
H.~Kesten.
\newblock The critical probability of bond percolation on the square lattice
  equals 1/2.
\newblock {\em Comm. Math. Phys.}, 74(1):41--59, 1980.

\bibitem{Kesten1982}
H.~Kesten.
\newblock {\em Percolation theory for mathematicians}.
\newblock Birkh\"{a}user, 1982.

\bibitem{Liggett1997}
T.~M. {Liggett}, R.~H. {Schonmann}, and A.~M. {Stacey}.
\newblock Domination by product measures.
\newblock {\em Ann. Prob.}, 25(1):71--95, 1997.

\bibitem{Meester1995}
R.~Meester.
\newblock Equality of critical densities in continuum percolation.
\newblock {\em J. Appl. Prob.}, 32(1):90--104, 1995.

\bibitem{Meester1996}
R.~Meester and R.~Roy.
\newblock {\em Continuum Percolation}.
\newblock Cambridge University Press, Cambridge, 1996.

\bibitem{Oskouyi2014}
A.~{Oskouyi}, U.~{Sundararaj}, and P.~{Mertiny}.
\newblock Effect of temperature on electrical resistivity of carbon nanotubes
  and graphene nanoplatelets nanocomposites.
\newblock {\em J. Nanotechnol. Eng. Med.}, 5(4):044501--044501, 2014.

\bibitem{Penrose2003}
M. D.~{Penrose}.
\newblock {\em Random geometric graphs}.
\newblock Oxford University Press, New York, 2003.

\bibitem{Penrose1991}
M. D.~{Penrose}.
\newblock On a continuum percolation model.
\newblock {\em Adv. Appl. Prob.}, 23(3):536--556, 1991.

\bibitem{Pike1974}
G.~E. Pike and C.~H. Seager.
\newblock Percolation and conductivity: A computer study. i.
\newblock {\em Phys. Rev. B}, 10(4):1421--1434, 1974.

\bibitem{RRoy1990}
R.~Roy.
\newblock The russo-seymour-welsh theorem and the equality of critical
  densities and the "dual"critical densities for continuum percolation on
  $\mathbb{R}^2$.
\newblock {\em Ann. Prob.}, 18(4):1563--1575, 1990.

\bibitem{RRoy1991}
R.~Roy.
\newblock Percolation of Poisson sticks on the plane.
\newblock {\em Probab. Th. Rel. Fields}, 89(4):503--517, 1991.

\bibitem{Russo1978}
L.~Russo.
\newblock A note on percolation.
\newblock {\em Z.Wahrsch. verw. Geb.}, 21:39--48, 1978.

\bibitem{SchneiderWeil2008}
R. ~{Schneider} and W.~{ Weil}.
\newblock {\em Stochastic and Integral Geometry}.
\newblock Springer, 2008.

\bibitem{Serre2013}
P.~{Serre}, C.~{Ternon}, V.~{Stambouli}, P.~{Periwal}, and T.~{Barona}.
\newblock Fabrication of silicon nanowire networks for biological sensing.
\newblock {\em Sensors and Actuators B: Chemical}, 182:390--395, 2013.

\bibitem{SW1978}
P.D.~{Seymour} and D. J. A.~{Welsh}.
\newblock Percolation probabilities on the square lattice.
\newblock {\em Ann. Discr. Math.}, 3:227--245, 1978.

\bibitem{Tassion2016}
V.~{Tassion}.
\newblock Crossing probabilities for Voronoi percolation.
\newblock {\em Ann. Prob}, 44(5):3385--3398, 2016.

\bibitem{Teixeira2017}
A.~Teixeira and D.~Ungaretti.
\newblock Ellipses percolation.
\newblock {\em J. Stat. Phy.}, 168(2):369--393, 2017.

\bibitem{Hofstad2017}
R. v. d. Hofstad.
\newblock {\em Random Graphs and Complex Networks}, volume~1 of {\em Cambridge
  Series in Statistical and Probabilistic Mathematics}.
\newblock Cambridge University Press, Cambridge, 2017.


\end{thebibliography}
\end{document}